\documentclass[11pt,a4paper,twoside]{amsart}
\usepackage{amsmath,amssymb,amsfonts,amsthm}
\usepackage{color,slashbox,textcomp,times,multirow}
\usepackage{framed,fancybox,empheq,lscape,graphics,footnote}
\usepackage[a4paper]{geometry}
\usepackage[colorlinks=true]{hyperref}

\usepackage{tikz}
\usetikzlibrary{calc}

\newtheorem{theorem}{Theorem}
\newtheorem{definition}[theorem]{Definition}

\newtheorem{proposition}[theorem]{Proposition}
\newtheorem{remark}[theorem]{Remark}
\newtheorem{example}{Example}
\newcommand{\mc}[3]{\multicolumn{#1}{#2}{#3}}

\title{Algebraic Multilevel Preconditioning in Isogeometric Analysis:\\
Construction and Numerical Studies
}

\author{K.P.S.~Gahalaut \and S.K.~Tomar \and J.K.~Kraus}
\address{Johann Radon Institute for Computational and Applied Mathematics,
Austrian Academy of Sciences \\ Altenbergerstrasse 69, 4040 Linz, Austria}
\email[Corresponding author]{krishan.gahalaut@ricam.oeaw.ac.at}

\date{June 21, 2013}


\begin{document}


\begin{abstract}
We present algebraic multilevel iteration (AMLI) methods for isogeometric discretization of scalar second order elliptic problems. The construction of coarse grid operators and hierarchical complementary operators are given. Moreover, for a uniform mesh on a unit interval, the explicit representation of B-spline basis functions for a fixed mesh size $h$ is given for $p=2,3,4$ and for  $C^{0}$- and $C^{p-1}$-continuity. The presented methods show $h$- and (almost) $p$-independent convergence rates. Supporting numerical results for convergence factor and iterations count for AMLI cycles ($V$-, linear $W$-, nonlinear $W$-) are provided. Numerical tests are performed, in two-dimensions on square domain and quarter annulus, and in three-dimensions on quarter thick ring.
\end{abstract}

\maketitle

\section{Introduction}
\label{sec:Intro}

The IsoGeometric Analysis (IGA), proposed by
Hughes et al. in \cite{HughesCB-05}, has received great attention
in the computational mechanics community. The concept has
the capability of leading to large steps forward in computational efficiency since
effectively, the process of re-meshing is either eliminated or greatly
suppressed. The geometry description of the underlying domain
is adopted from a Computer Aided Design (CAD) parametrization which is
usually based on Non-Uniform Rational B-splines (NURBS), and the same basis functions are 
employed to approximate the physical solution. Since its introduction, 
IGA techniques have been studied and applied in diverse fields, see e.g., \cite{AuricchioBHRS-10,BazilevsCHZ-08,BazilevsCCEHLSS-10,BuffaSV-10,BuffaRSV-11,CottrellRBH-06,GomezHNC-10,HughesRS-10}.
Moreover, some theoretical aspects such as approximation properties and condition number estimates have been studied, see 
\cite{BazilevsBCHS-06,BeiroBRS-11,BeiroCS-12,GahalautT-12}.
The isogeometric methods, depending on various choices of basis functions, have
shown several advantages over standard Finite Element Methods (FEM).
For example, some common geometries arising in engineering and applied sciences, 
such as circles or ellipses, are  exactly represented, 
and complicated geometries are represented more accurately than traditional polynomial based approaches.
When we compare NURBS based isogeometric analysis with standard 
Lagrange polynomials based finite element analysis, 
it leads to qualitatively more accurate results \cite{CHughesB-09}. 
Another limitation of finite element analysis is that it suits well for $C^0$ continuous interpolation, but
for $C^1$ or higher order interpolation finite elements are complicated 
and expensive to construct. IGA offers $C^{p-k}$-continuous interpolation for $p$-degree 
basis functions with knot multiplicity $k$. Moreover, the ease in building spaces with high
inter-element regularity allows for rather small problem sizes (in terms
of degrees of freedom) with respect to standard FEM with the same approximation properties. This implies
that, in general, for same approximation properties IGA stiffness and mass matrices are smaller than
the corresponding finite element ones. However, isogeometric matrices are denser than the FEM matrices in realistic problems of interest,
and their condition numbers grow quickly with the inverse of mesh
size $h$ and the polynomial degree $p$.  A detailed study of condition number estimates for the
stiffness matrix and mass matrix arising in isogeometric discretizations is given in \cite{GahalautT-12}. 
For the $h$-refinement, the condition number of the stiffness matrix is bounded from above and below by a constant times $ h^{-2}$,
and the condition number of the mass matrix is uniformly bounded.
For the $p$-refinement, the condition number is bounded above by
$p^{2d}4^{pd}$ and $p^{2(d-1)}4^{pd}$ for the stiffness matrix and the mass matrix, respectively, where $d$ is the dimension of the problem. 
As a consequence, the cost of solving the linear system of equations arising from the isogeometric discretization, particularly using iterative solvers,
becomes an important issue. Therefore, there is currently a growing
interest in the design of efficient preconditioners for IGA discrete
problems, in both the mathematical and the engineering
communities. Multigrid methods for IGA have been introduced for two and three dimensional elliptic problems 
by the authors in \cite{GahalautKT-12}, and tearing and interconnecting methods for isogeometric analysis are 
discussed in \cite{KleissPJT-12}. Other recent work on solvers for IGA studied overlapping additive Schwarz methods \cite{BeiroCPS-12,BeiroCPS-12-2} and balancing
domain decomposition by constraints methods \cite{BeiroCPS-12-1}. 
Some issues arising in using direct solvers have been investigated in \cite{CollierPDPC-12}.
The results, we presented in \cite{GahalautKT-12}, show optimal convergence rate with respect to the mesh parameter $h$. 
However, for discretizations based on higher degree polynomials,  the convergence rate are quickly deteriorated.
In this paper we discuss the construction of linear solvers which provide not only $h$-independent convergence rates but also exhibit (almost) independence on $p$. The presented optimal order solvers are based on algebraic multilevel iteration (AMLI) methods.

AMLI methods were introduced by Axelsson and Vassilevski
in a series of papers \cite{AxelssonV-89,AxelssonV-90,AxelssonV-91,AxelssonV-94}. The AMLI methods, which are recursive extensions of two-level
multigrid methods for FEM \cite{AxelssonG-83}, have been extensively analyzed in the context of conforming and nonconforming
FEM (including discontinuous Galerkin methods). For a detailed
systematic exposition of AMLI methods, see the monographs \cite{KrausBook,Vassilevski-08}. To reduce the overall complexity
of AMLI methods (to achieve optimal computational complexity), various stabilization techniques can be
used. In the original work \cite{AxelssonV-89,AxelssonV-90}, the stabilization was achieved by employing properly shifted and scaled
Chebyshev polynomials. This approach requires the computation of polynomial coefficients which depends
on the bounds of the eigenvalues of the preconditioned system. Alternatively, some inner iterations
at coarse levels can be used to stabilize the outer iterations, which lead to parameter-free
AMLI methods \cite{AxelssonV-91,AxelssonV-94,Kraus-02,Notay-02}. These methods utilize a sequence of coarse-grid problems that
are obtained from repeated application of a natural (and simple) hierarchical basis transformation, which
is computationally advantageous. Moreover, the underlying technique of these methods often requires
only a few minor adjustments (mainly two-level hierarchical basis transformation) even if the underlying
problem changes significantly.

In this article we consider the scalar second order elliptic equation as our model problem. 
Let $\Omega \subset \mathbb{R}^{d}, d = 2,3$, be an open, bounded and connected Lipschitz domain with Dirichlet boundary $\partial \Omega$. We consider
\begin{equation}
\label{eq:Poisson}
-\nabla \cdot (\mathcal{A} \nabla u) = f \quad \mathrm{in~} \Omega, \qquad u = u^D \quad \mathrm{on~} \partial \Omega,
\end{equation}
where $\mathcal{A}(x)$ is a uniformly bounded function for $x \in \Omega$. Let $V^0 \subset H^{1}(\Omega)$ denote the space of test functions which vanish on $\partial \Omega$, and $V^D = V^0 + u^D \subset H^{1}(\Omega)$ denote the set which contains the functions fulfilling the Dirichlet boundary condition on $\partial \Omega$. By $V^0_h \subset V^0$ and $V^D_h \subset V^D$ we denote the finite-dimensional spaces of the B-spline (NURBS) basis functions.
Introducing the bilinear form $a(\cdot , \cdot)$ and the linear form $f(\cdot)$ as
\begin{equation}
\label{eq:bilin_lin}
a(u,v) = \int_{\Omega} \mathcal{A} \nabla u \cdot \nabla v~dx, \quad
f(v) = \int_{\Omega} f~v~dx,
\end{equation} 
the Galerkin formulation of this problem reads:\\
Find $u_h \in V^D_h$ such that
\begin{equation}
\label{eq:VarProb_h}
a(u_h,v_h) = f(v_h) \quad \mathrm{for~all~} v_h \in V^0_h.
\end{equation}
It is well known that \eqref{eq:VarProb_h} is a well-posed problem and has a unique solution.
%
By approximating $u_h$ and $v_h$ using B-splines (NURBS) the variational formulation \eqref{eq:VarProb_h} is transformed in to a set of linear algebraic equations
\begin{equation}
\label{eq:LinSys}
A \bf u=f,
\end{equation}
where $A$ denotes the stiffness matrix obtained from the bilinear form $a(\cdot,\cdot)$, i.e.
 \begin{equation*}
A = (a_{i,j}) = (a(N_i,N_j)), \quad i,j=1,2,3,....,n_h,
\end{equation*}
$\bf u$ denotes the vector of unknown degrees of freedom (DOF), and $\bf f$ denotes the right hand side (RHS) vector from the known data of the problem.
Clearly, $A$ is a real symmetric positive definite matrix. 

The rest of the paper is organized as follows. In Section 2 we briefly review
the basics of B-splines and NURBS. An explicit representation of basis functions is also given in this section. 
The description of multilevel representation of B-splines (NURBS) is given in Section 3.
A brief description of AMLI methods is given in Section 4.
We then construct the isogeometric hierarchical spaces in Section 5. Numerical study of space splitting techniques is discussed in Section 6.
The results of AMLI methods for several numerical experiments in two- and three-dimensions are presented in Section 7. Finally, some conclusions are drawn in Section 8.

\section{B-splines and NURBS}
\label{sec:BS_NURBS}

\subsection{B-splines}
We first recall the definition of B-splines, see e.g. \cite{Boor-78,PieglTiller-97, Rozers-01,Schumaker-07}.
\begin{definition}
\label{Def:BSplines}
Let $\Xi_{1}=\{\xi_i : i = 1, ..., n + p+1\}$ be a non-decreasing sequence of real numbers, called the $knot$ $vector$, where $\xi_i$ is the $i^{th}$ knot, 
$p$ is the polynomial degree, and $n$ is the number of basis function.
With a knot vector in hand, the B-spline basis functions, denoted by $ N^p_{i}(\xi)$, are (recursively) defined starting with a piecewise constant 
\begin{subequations}
\label{eq:DefBSplines}
\begin{alignat}{3}
B^0_{i}(\xi) &= \begin{cases} 1 & \text{if $\xi \in [\xi_i,\xi_{i+1})$,} \\ 0 &\text{otherwise,} \end{cases}\\
\quad B^p_{i}(\xi) & =\frac{\xi-\xi_i}{\xi_{i+p}-\xi_i} B^{p-1}_{i}(\xi) +\frac{\xi_{i+p+1}-\xi}{\xi_{i+p+1}-\xi_{i+1}} B^{p-1}_{i+1}(\xi),
\end{alignat}
\end{subequations}
where $ 0 \leq i \leq n, p \geq 1$, and $\displaystyle \frac{0}{0}$ is considered as zero.
\end{definition}
The above expression is usually referred as the {\it{Cox-de Boor recursion formula}}, see e.g. \cite{Boor-78}.
For a B-spline basis function of degree $p$, an interior knot can be repeated at most $p$ times, and the boundary knots can be repeated at most $p+1$ times.
A knot vector for which the two boundary knots are repeated $p+1$ times is said to be open.
In this case, the basis functions are interpolatory at the first and the last knot.
Important properties of the B-spline basis functions include nonnegativity, partition of unity, local support and $C^{p-k}$-continuity.
\begin{definition}
\label{Def:BSplineCurve}
A B-spline curve $C(\xi)$, is defined by
\begin{equation}
C(\xi)=\sum^{n}_{i=1} P_{i} B^p_{i}(\xi)
\end{equation}
where $\{P_i : i = 1, ..., n\}$ are the control points and $B^p_{i}$ are B-spline basis functions defined in (\ref{eq:DefBSplines}).
\end{definition}
The previous definitions are easily generalized to the higher dimensional cases by means of tensor product.
Using tensor product of one-dimensional B-spline functions, a B-spline surface $S(\xi,\eta)$ is defined as follows:
\begin{equation}
\label{eq:BSurface}
S(\xi,\eta) = \sum_{i=1}^{n_{1}} \sum_{j=1}^{n_{2}} B_{i,j}^{p_{1},p_{2}} (\xi,\eta) P_{i,j},
\end{equation}
where $P_{i,j}$, $i = 1,2, \ldots, n_{1}$, $j = 1,2, \ldots, n_{2}$, denote the control points,
$B_{i,j}^{p_{1},p_{2}}$ is the tensor product of B-spline basis functions $B_{i}^{p_{1}}$
and $B_{j}^{p_{2}}$, and $\Xi_{1} = \{\xi_{1}, \xi_{2}, \ldots, \xi_{n_{1}+p_{1}+1}\}$
and $\Xi_{2} = \{\eta_{1}, \eta_{2}, \ldots, \eta_{n_{2}+p_{2}+1}\}$ are the corresponding knot vectors.
Similarly, B-spline solids can be defined by a three-dimensional tensor product.

\subsection{NURBS}

While B-splines (polynomials) are flexible and have many nice properties for curve design, they are also incapable of exactly representing curves such as circles, ellipses, etc.
Such limitations are overcome by NURBS functions.
Rational representation of conics originates from projective geometry and requires additional parameters called weights, which we shall denote by $w$. 
Let $\{P^{w}_{i}\}$ be a set of control points for a projective B-spline curve in $\mathbb{R}^{3}$.
For the desired NURBS curve in $\mathbb{R}^{2}$, the weights and the control points are derived by the relations
\begin{equation}
\label{eq:NWts}
w_{i} = (P^{w}_{i})_{3}, \qquad (P_{i})_{d} = (P^{w}_{i})/w_{i}, \quad d = 1,2,
\end{equation}
where $w_{i}$ is called the $i^{\mathrm{th}}$ weight and $(P_{i})_{d}$ is the $d^{\mathrm{th}}$-dimension
component of the vector $P_{i}$. 
The weight function $w(\xi)$ is defined as
\begin{equation}
\label{eq:NWFunc}
w(\xi) = \sum_{i = 1}^{n}B_i^p(\xi) w_{i}.
\end{equation}
Then, the NURBS basis functions and curve are defined by
\begin{equation}
\label{eq:NCurve}
N_i^p(\xi)  = \frac{B_i^p(\xi)  w_{i}}{w(\xi)}, \qquad C(\xi) = \sum_{i=1}^{n} N_i^p(\xi)  P_{i}.
\end{equation}
The NURBS surfaces are analogously defined as follows
\begin{equation}
\label{eq:NSurface}
S(\xi,\eta) = \sum_{i=1}^{n_{1}} \sum_{j=1}^{n_{2}} N_{i,j}^{p_{1}, p_{2}} (\xi,\eta) P_{i,j},
\end{equation}
where $N_{i,j}^{p_{1},p_{2}}$ is the tensor product of NURBS basis functions $N_{i}^{p_{1}}$ and $N_{j}^{p_{2}}$.
Similarly, NURBS solids can be defined by a three-dimensional tensor product.
NURBS functions also satisfy the properties of B-spline functions.
For a detailed exposition see, e.g. \cite{Boor-78, PieglTiller-97, Rozers-01, Schumaker-07}.

\subsection{Explicit Representation for B-splines}

The recursive form of B-spline basis functions, given by \eqref{eq:DefBSplines}, is elegant and concise, and is presented in all the IGA related references, see e.g., \cite{Boor-78, PieglTiller-97, Rozers-01, Schumaker-07}. However, this form may not be the most efficient from computational point of view, specially when dealing with large knot vectors. 
Therefore, in the CAD community, there have been considerable efforts for efficient NURBS evaluation techniques, for latest advances see, e.g., optimized GPU evaluation of NURBS curves and surfaces \cite{Krishnamurti-09} and references therein.
To the best of authors' knowledge, within the IGA literature there is no reference on the explicit representation of B-splines for a given mesh size $h$. However, there are situations, e.g. in academic problems, where having an explicit representation of B-spline basis functions is of significant importance. Therefore, we present the explicit form of B-splines in terms of the mesh size (knot-span) $h$. Having an explicit form of basis functions is also advantageous in devising inter-grid transfer operators for multigrid and multilevel iterative solvers.  For brevity reasons, we restrict ourselves to a unit interval with equal spacing. Moreover, as most of the NURBS based designs in engineering use polynomial degree $p=2$ and $3$, we will confine ourselves up to $p=4$ with $C^0$ and $C^{p-1}$ continuous basis functions. 

\subsubsection{$C^{p-1}$-continuity}
\label{subsec:ExBp-1} 

We first consider the $C^{p-1}$ continuous case as this is the default case for knot vector with non-repeated internal knots. For B-spline functions with $p = 0$ and $p = 1$, we have the same representation as for standard piecewise
constant and linear finite element functions, respectively. 
Quadratic B-spline basis functions,
however, differ from their FEM counterparts. They are each identical but shifted related to
each other, whereas the shape of a quadratic finite element function depends on whether it
corresponds to an internal node or an end node. This ``homogeneous'' pattern continues for the
B-splines with higher-degrees. 

We are interested to give an explicit representation for uniform B-spline basis functions defined on a knot vector $E_k$ at any given level $k$, where $k = 1,2,...,L,$ with spacing $h$ $(=1/n)$, where $n$ is the total number of knot spans. We shall use the notation $B^{p,r}_{k,i}$ for B-splines, where superscripts represent the polynomial degree and the regularity of basis  functions, respectively, and the subscripts represent the level and the number of basis function, respectively. We start with level $1$ with only one element. Using the definition from \eqref{eq:DefBSplines}, at level $1$ the B-spline basis functions of degree $p=2$ on the knot vector $E_1 = \{0, 0, 0, 1, 1, 1\}$ are defined as follows 
\begin{equation}
\label{eq:bsplinep2E1}
\begin{split}
B^{2,{p-1}}_{1,1} & = (1-x)^2, \quad 0 \le x \le1,\\ 
B^{2,{p-1}}_{1,2} & = 2x(1-x), \quad 0 \le x \le1,\\ 
B^{2,{p-1}}_{1,3} & = x^2,  \quad  0 \le x \le1.
\end{split}
\end{equation}
The mesh refinement takes place by inserting the knots. We consider uniform refinement of $E_1$, i.e. inserting knots at the mid point of the knot values. At the next level $k=2$, the basis functions on refined knot vector  $E_2 = \{0, 0, 0, \frac{1}{2},1, 1, 1\}$ are given by
\begin{equation*}
\begin{split}
B^{2,{p-1}}_{2,1}& =
\begin{cases}     
    (1-2x)^2,  \quad  0 \le x < \frac{1}{2},\\
    0,  \quad  \frac{1}{2} \le x \le 1,
\end{cases} \\ 
B^{2,{p-1}}_{2,2}&=
\begin{cases}     
   2x(2-3x), \quad  0 \le x < \frac{1}{2},\\
     2(1-x)^2,\quad   \frac{1}{2} \le x \le 1,
\end{cases} \\ 
B^{2,{p-1}}_{2,3}&=
\begin{cases}     
   2x^2,  \quad   0 \le x < \frac{1}{2},\\
     -2+8x-6x^2,  \quad  \frac{1}{2}\le x  \le  1,
\end{cases}
\end{split}
\end{equation*}
\begin{equation}
\begin{split}
B^{2,{p-1}}_{2,4}&=
\begin{cases}     
   0,   \quad 0 \le x < \frac{1}{2},\\
     (1-2x)^2,   \quad  \frac{1}{2} \le x \le1.
\end{cases} 
\end{split}
\end{equation}
Further refinements take place in a similar way, i.e., starting with $E_1$, a single knot span, in the knot span $E_k$ we will thus have $2^{k-1}$ knot spans.  The explicit representation of B-splines at level $k$, where $k \ge 3,$ is given by 
\begin{equation*}
\begin{split}
%
B^{2, {p-1}}_{k,1} & = \frac{1}{h^2}(h-x)^2, \quad 0 \le x < h,\quad h \le 1,\\ 
B^{2, {p-1}}_{k,2} & =
\begin{array}{c|c}
\begin{cases}  \displaystyle 
    \frac{1}{2h^2}x(4h-3x),   & \quad 0 \le x < h, \vspace{1mm}\\ \displaystyle
    \frac{1}{2h^2}(2h-x)^2,  & \quad h \le x<2h,
\end{cases} &  \displaystyle \text{ for }h \le \frac{1}{2}, 
\end{array} 
\end{split}
\end{equation*}
\begin{equation}
\begin{split}
B^{2, {p-1}}_{k,3+i} & =
\begin{cases}  \displaystyle 
    \frac{1}{2h^2}(x-ih)^2,   & \quad ih \le x < (i+1)h, \vspace{1mm}\\ \displaystyle
     \frac{-3}{2}+\frac{3}{h}(x-ih)-\frac{1}{h^2}(x-ih)^2,  & \quad (i+1)h \le x < (i+2)h,  \vspace{1mm}\\ \displaystyle
 \frac{1}{2h^2}(3h-(x-ih))^2,  & \quad (i+2)h \le x<(i+3)h,
\end{cases}\\
&\hspace{3.5cm}\mathrm{ where} \quad i =0, 1,2,3,...,(1/h)-3, \text{ and } h \le 1/4.\\
B^{2, {p-1}}_{k,n+p-1} & =
\begin{array}{c|c}
\begin{cases}  \displaystyle 
    \frac{1}{2h^2}(-1+2h+x)^2,    \quad 1-2h \le x < 1-h, \vspace{1mm}\\ \displaystyle
    \frac{-1}{2h^2}(3-4h+2(2h-3)x+3x^2),   \quad 1-h \le x \le 1,
\end{cases}&  \displaystyle \text{ for }h \le \frac{1}{2}, 
\end{array} \\
B^{2, {p-1}}_{k,n+p} & = \frac{1}{h^2}(h-(1-x))^2, \quad 1-h \le x \le 1, \quad h \le 1.
\end{split}
\end{equation}
For higher degree polynomials, we can define the explicit representation in a similar way. Again using the definition \eqref{eq:DefBSplines} of B-splines, for $p=3$, at first level $k=1$, the basis functions with $C^{p-1}$-continuity are given as follows
\begin{equation}
\label{eq:bsplinep3E1}
\begin{split}
B^{3,{p-1}}_{1,1} & = (1-x)^3,  \quad 0 \le x \le1,\\
B^{3, {p-1}}_{1,2} & = 3x(1-x)^2,  \quad 0 \le x \le1,\\
B^{3, {p-1}}_{1,3} & = 3x^2(1-x),  \quad 0 \le x \le1, \\ 
B^{3, {p-1}}_{1,4} & = x^3, \quad  0 \le x \le1.
\end{split}
\end{equation}
At level $2$, we have the following basis functions
\begin{equation}
\begin{split}
B^{3, {p-1}}_{2,1}& =
\begin{cases}     
    (1-2x)^3,   \quad 0 \le x < \frac{1}{2},\\
    0,   \quad \frac{1}{2} \le x \le 1,
\end{cases} \\ 
B^{3, {p-1}}_{2,2}&=
\begin{cases}     
    2x(3-9x+7x^2) , \quad  0 \le x < \frac{1}{2},\\
 2(1-x)^3  , \quad  \frac{1}{2} \le x \le 1,
\end{cases} \\
B^{3, {p-1}}_{2,3}&=
\begin{cases}     
  2x^2(3-4x),   \quad  0 \le x < \frac{1}{2},\\
      2(-1+x)^2(-1+4x),  \quad  \frac{1}{2}\le x  \le  1,
\end{cases}\\
B^{3, {p-1}}_{2,4}&=
\begin{cases}     
   2x^3,   \quad  0 \le x < \frac{1}{2},\\
       2-12x+24x^2-14x^3, \quad  \frac{1}{2} \le x \le 1,
\end{cases} \\ 
B^{3, {p-1}}_{2,5}&=
\begin{cases}     
   0,   \quad  0 \le x < \frac{1}{2},\\
     (-1+2x)^3,   \quad  \frac{1}{2} \le x \le 1.
\end{cases} 
\end{split}
\end{equation}
For all other levels $k$, where $k\ge3,$ the basis functions are defined below
\begin{equation*}
\begin{split}
%
B^{3, {p-1}}_{k,1} & = \frac{1}{h^3}(h-x)^3, \quad 0 \le x < h,\quad h \le 1,\\
B^{3, {p-1}}_{k,2} & =
\begin{array}{c|c}
\begin{cases}  \displaystyle 
    \frac{x}{h}\left(3- \frac{9}{2}\frac{x}{h}+\frac{7}{4}\frac{x^2}{h^2}\right),   & \quad 0 \le x < h, \vspace{1mm}\\ \displaystyle
    \frac{1}{4h^3}(-2h+x)^3,  & \quad h \le x<2h, 
\end{cases}  & \displaystyle\text{ for }h \le \frac{1}{2}, 
\end{array}  \\
%
%
B^{3, {p-1}}_{k,3} & =
\begin{array}{c|c}
\begin{cases}  \displaystyle 
  \frac{1}{6}  \frac{x^2}{h^2}\left(9-\frac{11}{2}\frac{x}{h}\right),   & \quad 0 \le x < h, \vspace{1mm}\\ \displaystyle
 \frac{-3}{2}+\frac{9}{2}\frac{x}{h} -3\frac{x^2}{h^2} +\frac{7}{4}\frac{x^3}{h^3},         & \quad h \le x<2h, \vspace{1mm}\\ \displaystyle
    \frac{1}{6h^3}(-3h+x)^3,  & \quad 2h \le x<3h, 
\end{cases} &  \displaystyle \text{ for }h \le \frac{1}{4}, 
\end{array} \\
%
%
%
%
B^{3, {p-1}}_{k,4+i} & =
\begin{cases}  \displaystyle 
\frac{1}{6h^3}(x-ih)^3,   \hfill ih \le x < (i+1)h, \vspace{1mm}\\ \displaystyle
\frac{2}{3}-\frac{2}{h}(x-ih)+\frac{1}{2h^2}(x-ih)^2 - \frac{1}{2h^3}(x-ih)^3,   
 \hfill (i+1) \le x<(i+2)h, \\ \displaystyle
\frac{-22}{3}+\frac{10}{h}(x-ih)-\frac{4}{h^2}(x-ih)^2+\frac{1}{2h^3}(x-ih)^3,  
 \hfill (i+2)h \le x<(i+3)h, \\ \displaystyle
\frac{32}{h} \left(1-\frac{(x-ih)}{4h}\right)^3,   \hfill  (i+3)h \le x<(i+4)h,
\end{cases}\\
&\hspace{4cm} \hfill \text{ where} \quad i =0, 1,2,3,...,(1/h)-4, \text{ and } h \le \frac{1}{4}, \\
B^{ 3, {p-1}}_{k,n+p-2} & = 
\begin{cases}  \displaystyle 
   \frac{1}{6h^3}(-3h+(1-x))^3 ,   &  1-3h \le x < 1-2h, \vspace{1mm}\\ \displaystyle
 \frac{-3}{2}+\frac{9}{2}\frac{(1-x)}{h} -3\frac{(1-x)^2}{h^2} +\frac{7}{4}\frac{(1-x)^3}{h^3},         & 1-2h \le x<1-h, \vspace{1mm}\\ \displaystyle
 \frac{1}{6}  \frac{(1-x)^2}{h^2}\left(9-\frac{11}{2}\frac{(1-x)}{h}\right) ,  &  1-h \le x \le 1, 
\end{cases}\\
& \hspace{10cm} \hfill \text{for }  h \le \frac{1}{4}, 
\end{split}
\end{equation*}
\begin{equation}
\begin{split}
B^{3, {p-1}}_{k,n+p-1} & =
\begin{array}{c|c}
\begin{cases}  \displaystyle 
   \frac{1}{4h^3}(-2h+(1-x))^3 , \quad \quad \quad \quad \quad \quad 1-2h \le x < 1-h, \vspace{1mm}\\ \displaystyle
  \frac{(1-x)}{h}\left(3- \frac{9}{2}\frac{(1-x)}{h}+\frac{7}{4}\frac{(1-x)^2}{h^2}\right) , \quad 1-h \le x<1, 
\end{cases} & \displaystyle\text{for }h \le \frac{1}{2}, 
\end{array}\\
%
%
%
B^{3,{p-1}}_{k,n+p} & = \frac{1}{h^3}(h-(1-x))^3, \quad 1-h \le x \le 1,\quad h \le 1.\\ 
\end{split}
\end{equation}
Finally, we give the explicit representation of basis functions for $p=4$ with $C^{p-1}$-continuity.
At level $1$, with knot $E_1$, the B-spline basis functions of degree $p=4$ are given by
\begin{equation}
\label{eq:bsplinep4E1}
\begin{split}
B^{4, {p-1}}_{1,1} & = (1-x)^4, \quad  0 \le x \le1,\\  
B^{4, {p-1}}_{1,2} & = 4x(1-x)^3, \quad  0 \le x \le1,\\ 
B^{4,{p-1}}_{1,3} & = 6x^2(1-x)^2, \quad  0 \le x \le1,\\
B^{4, {p-1}}_{1,4} & = 4x^3(1-x),  \quad 0 \le x \le1, \\  
B^{4, {p-1}}_{1,5} & = x^4, \quad  0 \le x \le1.
\end{split}
\end{equation}
The B-splines on second level $k=2$ with knot $E_2$ are defined as follows
\begin{equation}
\begin{split}
B^{4, {p-1}}_{2,1}& =
\begin{cases}     
    (1-2x)^4,   \quad 0 \le x < \frac{1}{2},\\
    0,   \quad \frac{1}{2} \le x \le 1,
\end{cases} \\ 
B^{4, {p-1}}_{2,2}&=
\begin{cases}     
    2x(4-18x+28x^2-15x^3) , \quad  0 \le x < \frac{1}{2},\\
 2(1-x)^4  , \quad  \frac{1}{2} \le x \le 1,
\end{cases} \\ 
B^{4, {p-1}}_{2,3}&=
\begin{cases}     
  2x^2(6-16x+11x^2),   \quad  0 \le x < \frac{1}{2},\\
      2(1-x)^3(-1+5x),  \quad  \frac{1}{2}\le x \le 1,
\end{cases}\\ 
B^{4, {p-1}}_{2,4}&=
\begin{cases}     
  2x^3(4-5x),   \quad  0 \le x < \frac{1}{2},\\
      2(1-x)^2(1-6x+11x^2),  \quad  \frac{1}{2}\le x \le 1,
\end{cases}\\ 
B^{4, {p-1}}_{2,5}&=
\begin{cases}     
   2x^4,   \quad  0 \le x < \frac{1}{2},\\
   -2+16x-48x^2+64x^3-30x^4, \quad  \frac{1}{2} \le x \le 1,
\end{cases} \\ 
B^{4, {p-1}}_{2,6}&=
\begin{cases}     
   0,   \quad  0 \le x < \frac{1}{2},\\
     (1-2x)^4,   \quad  \frac{1}{2} \le x \le 1.
\end{cases} 
\end{split}
\end{equation}
At all other levels $k$, where $k\ge3,$ the basis functions of degree $p=4$ with $C^{p-1}$-continuity are given by
\begin{equation*}
\begin{split}
%
B^{4, {p-1}}_{k,1} & = \frac{1}{h^4}(h-x)^4, \quad 0 \le x < h,\quad h \le 1,\\ 
B^{4, {p-1}}_{k,2} & =
\begin{array}{c|c}
\begin{cases}  \displaystyle 
    \frac{-4x}{h}\left(-1+ \frac{9}{4}\frac{x}{h}-\frac{7}{4}\frac{x^2}{h^2}+\frac{15}{32}\frac{x^3}{h^3}\right),   & \quad 0 \le x < h, \vspace{1mm}\\ \displaystyle
    \frac{1}{8h^4}(2h-x)^4,  & \quad h \le x<2h, 
\end{cases} & \displaystyle \text{ for } h \le \frac{1}{2}, 
\end{array} \\
B^{4, {p-1}}_{k,3} & =
\begin{array}{c|c}
\begin{cases}  \displaystyle 
  \frac{1}{9}  \frac{x^2}{h^2}\left(27-33\frac{x}{h}+\frac{85}{8}\frac{x^2}{h^2}\right),   & \quad 0 \le x < h, \vspace{1mm}\\ \displaystyle
 \frac{-3}{2}+6\frac{x}{h} -6\frac{x^2}{h^2} +\frac{7}{3}\frac{x^3}{h^3}-\frac{23}{72}\frac{x^4}{h^4},         & \quad h \le x<2h, \vspace{1mm}\\ \displaystyle
    \frac{1}{18h^4}(3h-x)^4,  & \quad 2h \le x<3h, 
\end{cases}& \displaystyle \text{ for } h \le \frac{1}{4}, 
\end{array}\\
%
%
%
B^{4, {p-1}}_{k,4} & =
\begin{array}{c|c}
\begin{cases}  \displaystyle 
  \frac{2}{3}\frac{x^3}{h^3}-\frac{25}{72}\frac{x^4}{h^4},  & \quad 0 \le x<h,\vspace{1mm} \\ \displaystyle
 \frac{2}{3}-\frac{8}{3}\frac{x}{h}+4\frac{x^2}{h^2} -2 \frac{x^3}{h^3}+\frac{23}{72}\frac{x^4}{h^4},         & \quad h \le x<2h,\vspace{1mm} \\ \displaystyle
 \frac{-22}{3}+\frac{40}{3}\frac{x}{h} -8\frac{x^2}{h^2} +2\frac{x^3}{h^3}-\frac{13}{72}\frac{x^4}{h^4},         & \quad 2h \le x<3h,\vspace{1mm} \\ \displaystyle
    \frac{1}{24h^4}(4h-x)^4,  & \quad 3h \le x<4h, 
\end{cases}& \displaystyle \text{ for } h \le \frac{1}{4}. 
\end{array}
%
%
\end{split}
\end{equation*}
\begin{equation*}
\begin{split}
%
%
%
B^{4,{p-1}}_{k,5+i} & =
\begin{cases}  \displaystyle 
    \frac{1}{24h^4}(x-ih)^4,  \hfill ih \le x<(i+1)h, \vspace{1mm} \\ \displaystyle
  \frac{1}{24}\left(-5+\frac{20}{h}(x-ih) -\frac{30}{h^2}(x-ih)^2 + \frac{20}{h^3}(x-ih)^3 - \frac{4}{h^4}(x-ih)^4 \right) , \\
  \hfill (i+1)h \le x<(i+2)h, \\ \displaystyle
   \frac{155}{24}-\frac{25}{2h}(x-ih) +\frac{35}{4h^2}(x-ih)^2 - \frac{5}{2h^3}(x-ih)^3 - \frac{1}{4h^4}(x-ih)^4  ,\\
   \hfill (i+2)h \le x<(i+3)h, \\ \displaystyle
   \frac{-655}{24}+\frac{65}{2h}(x-ih)-\frac{55}{4h^2}(x-ih)^2 +\frac{5}{2h^3}(x-ih)^3 - \frac{1}{6h^4}(x-ih)^4  , \\
  \hfill (i+3)h \le x<(i+4)h, \\ \displaystyle
    \frac{1}{24h^4}(5h-(x-ih))^4,   \hfill (i+4)h \le x<(i+5)h, 
\end{cases}\\
&\hspace{4cm}\mathrm{ where} \quad i =0, 1,2,3,...,(1/h)-5, \text{ and } h \le \frac{1}{8}, (k \ge 4),\\
%
%
B^{4, {p-1}}_{k,n+p-3} & =
\begin{cases}  \displaystyle 
    \frac{1}{24h^4}\left(4h-(1-x)\right)^4,   \hfill 1-4h \le x<1-3h, \vspace{1mm} \\ \displaystyle
 \frac{-22}{3}+\frac{40}{3}\frac{(1-x)}{h} -8\frac{(1-x)^2}{h^2} +2\frac{(1-x)^3}{h^3}-\frac{13}{72}\frac{(1-x)^4}{h^4},     \\
     \hfill 1-3h \le x<1-2h, \\ \displaystyle
 \frac{2}{3}-\frac{8}{3}\frac{(1-x)}{h}+4\frac{(1-x)^2}{h^2} -2 \frac{(1-x)^3}{h^3}+\frac{23}{72}\frac{(1-x)^4}{h^4},  \\
 \hfill 1-2h \le x<1-h, \\ \displaystyle
  \frac{2}{3}\frac{(1-x)^3}{h^3}-\frac{25}{72}\frac{(1-x)^4}{h^4},  \hfill 1-h \le x\le1, 
\end{cases}\\
&\hspace{9cm}\text{ for }  h \le \frac{1}{4}, \\
B^{4, {p-1}}_{k,n+p-2} & =
\begin{cases}  \displaystyle 
    \frac{1}{18h^4}(3h-(1-x))^4,  \hfill 1-3h \le x<1-2h, \vspace{1mm} \\\displaystyle
\frac{-3}{2}+6\frac{(1-x)}{h} -6\frac{(1-x)^2}{h^2} +\frac{7}{3}\frac{(1-x)^3}{h^3}-\frac{23}{72}\frac{(1-x)^4}{h^4}, \\
 \hfill 1-2h \le x<1-h, \\\displaystyle
  \frac{1}{9}  \frac{(1-x)^2}{h^2}\left(27-33\frac{(1-x)}{h}+\frac{85}{8}\frac{(1-x)^2}{h^2}\right),    \hfill 1-h \le x \le 1 ,
\end{cases}\\
&\hspace{9cm}\text{ for }  h \le \frac{1}{4}, 
\end{split}
\end{equation*}
\begin{equation}
\begin{split}
B^{4, {p-1}}_{k,n+p-1} & =
\begin{cases}  \displaystyle 
\frac{1}{8h^4}(2h-(1-x))^4,  \hfill 1-2h \le x<1-h, \vspace{1mm} \\\displaystyle
    \frac{-4(1-x)}{h}\left(-1+ \frac{9}{4}\frac{(1-x)}{h}-\frac{7}{4}\frac{(1-x)^2}{h^2}+\frac{15}{32}\frac{(1-x)^3}{h^3}\right), \\
  \hfill 1-h \le x \le 1, 
\end{cases}\\
&\hspace{9cm}\text{ for }  h \le \frac{1}{2}, \\
B^{4, {p-1}}_{k,n+p} & = \frac{1}{h^4}(h-(1-x))^4, \quad 1-h \le x \le 1, \quad h \le 1.
\end{split}
\end{equation}
Note that the expression $B^{4,{p-1}}_{k,5+i}$ is valid only for $k \ge 4$.

\subsection{$C^{0}$-continuity}
\label{subsec:ExB0}

To reduce the continuity of the basis functions across  element boundaries, the knot values are repeated upto a desired level. By repeating the internal knots $k$ times we get the $C^{p-k}$ continuous basis functions.
In the previous section we have given the explicit representation for $C^{p-1}$ continuity, which is the highest continuity for polynomial degree $p$. We now consider another extreme case, the lowest continuity, i.e. $C^0$ continuous basis functions.
At first level $k=1$ the $C^0$ continuous B-spline basis functions of degree $p=2,3,4$ on a knot $E_1 = \{0, 0, 0, 1, 1, 1\}$ 
are same as those of $C^{p-1}$ continuous B-spline basis functions of same degree, see \eqref{eq:bsplinep2E1}, \eqref{eq:bsplinep3E1}, and \eqref{eq:bsplinep4E1} respectively.

The explicit representation for $C^0$ continuous B-spline basis functions of degree $p=2$ at level $k$, where $k\ge2$ is given by 
\begin{equation}
\begin{split}
\label{eq:Bsplinesp2E1} 
B^{2, {0}}_{k,1} & = \frac{1}{h^2}(h-x)^2, \quad 0 \le x < h,\\ 
B^{2, {0}}_{k,2+2i} & = \frac{-2}{h^2}(x-ih)(h+(x-ih)),  \quad (i-1)h \le x < ih,\\ 
&\hspace{3cm}\text{ where } i = 0,1,2,3,...,(1/h)-1,\\
B^{2, {0}}_{k,3+2i} & =
\begin{cases}  \displaystyle 
    \frac{1}{h^2}(h+(x-ih))^2,   & \quad(i-1)h \le x < ih, \vspace{1mm}\\ \displaystyle
    \frac{1}{h^2}(-h+(x-ih))^2,  & \quad ih \le x<(i+1)h,
\end{cases}\\
&\hspace{3cm}\text{ where }  i = 0,1,2,3,...,\left((1/h)-2\right),\\
B^{2, {0}}_{k,np+1} & = \frac{1}{h^2}(h-(1-x))^2, \quad 1-h \le x \le 1.
\end{split}
\end{equation}
For $p=3$, the explicit representation for B-spline basis functions with $C^0$-continuity, at level $k$, where $k\ge2$  is given by
\begin{equation}
\begin{split}
\label{eq:Bsplinesp3E1} 
B^{3, {0}}_{k,1} & = \frac{1}{h^3}(h-x)^3, \quad 0 \le x < h,\\ 
B^{3, {0}}_{k,2+3i} & = \frac{3}{h}\left(-1+\frac{1}{h}(x-ih)\right)^2,  \quad ih \le x < (i+1)h,\\ 
&\hspace{3cm}\text{ where }  i = 0,1,2,3,...,(1/h)-1,\\
B^{3, {0}}_{k,3+3i} & = \frac{3}{h^2}(x-ih)^2 \left(1- \frac{1}{h}(x-ih)\right),  \quad ih \le x < (i+1)h,\\ 
&\hspace{3cm}\text{ where }  i = 0,1,2,3,...,(1/h)-1,\\
B^{3, {0}}_{k,4+3i} & =
\begin{cases}  \displaystyle 
    \frac{1}{h^3}(x-ih)^3,   & \quad ih \le x < (i+1)h, \vspace{1mm}\\ \displaystyle
  8\left(1-  \frac{1}{2h}(x-ih)\right)^3,  & \quad (i+1)h \le x<(i+2)h,
\end{cases}\\
&\hspace{3cm}\text{ where }  i = 0,1,2,3,...,\left((1/h)-2\right),\\
B^{3, {0}}_{k,np+1} & = \frac{1}{h^3}(h-(1-x))^3, \quad 1-h \le x \le 1.
\end{split}
\end{equation}
Finally, the explicit representation for $C^0$ continuous basis functions of degree $p=4$  at level $k$, where $k\ge2$ is given below
\begin{equation*}
\begin{split}
\label{eq:Bsplinesp4E1} 
B^{4, {0}}_{k,1} & = \frac{1}{h^4}(h-x)^4, \quad 0 \le x < h,\\ 
B^{4, {0}}_{k,2+4i} & = \frac{4}{h}(x-ih)\left(1-\frac{(x-ih)}{h}\right)^3,  \quad ih \le x < (i+1)h,\\ 
&\hspace{3cm}\text{ where } i = 0,1,2,3,...,(1/h)-1,\\
\end{split}
\end{equation*}
\begin{equation}
\begin{split}
B^{4, {0}}_{k,3+4i} & = \frac{6}{h^2}(x-ih)^2\left(1-\frac{(x-ih)}{h}\right)^2,  \quad ih \le x < (i+1)h,\\ 
&\hspace{3cm}\text{ where }  i = 0,1,2,3,...,(1/h)-1,\\
B^{4, {0}}_{k,4+4i} & = \frac{4}{h^3}(x-ih)^3\left(1-\frac{(x-ih)}{h}\right),  \quad ih \le x < (i+1)h,\\ 
&\hspace{3cm}\text{ where }  i = 0,1,2,3,...,(1/h)-1,\\
B^{4, {0}}_{k,5+4i} & =
\begin{cases}  \displaystyle 
    \frac{1}{h^4}(x-ih)^4,   & \quad ih \le x < (i+1)h, \vspace{1mm}\\ \displaystyle
  16\left(1-  \frac{1}{2h}(x-ih)\right)^4,  & \quad (i+1)h \le x<(i+2)h,
\end{cases}\\
&\hspace{3cm}\text{ where }   i = 0,1,2,3,...,\left((1/h)-2\right),\\
B^{4, {0}}_{k,np+1} & = \frac{1}{h^4}(h-(1-x))^4, \quad 1-h \le x \le 1.
\end{split}
\end{equation}

\section{Multilevel Representation of B-splines and NURBS}
\label{sec:MulSplines}

\subsection{Multilevel B-splines}

In this section, we study the multilevel structure of B-splines and NURBS spaces. This will be used in the construction of corresponding hierarchical spaces (i.e. splitting the fine space into coarse space and its hierarchical complement) in Section \ref{sec:hier_space}. 
For standard FEM, multilevel $h$ and $p$ representations can be built with nice approximation properties, see e.g., \cite{Carey-97}. 
The recent work on refinement strategies of B-splines can be found in \cite{SchillingerDS-11,SchillingerR-11} where local refinement techniques have been discussed. We focus on generating an explicit matrix form of transfer operators. 
For a two level setting, let $\mathcal{B}^{p,r}_{k-1}$ and $\mathcal{B}^{p,r}_{k}$ denote the B-spline spaces at coarse and fine level, respectively. Let $\{B^{p,r}_{k-1,i}, i = 1,2,...,n_{k-1}\}$ and $\{B^{p,r}_{k,i}, i = 1,2,...,n_k\}$ be the set of basis functions for coarse and fine space, respectively, i.e.
\begin{equation*}
\mathcal{B}^{p,r}_{k-1} = \mathrm{span}\{B^{p,r}_{k-1,1},B^{p,r}_{k-1,2},B^{p,r}_{k-1,3},...,B^{p,r}_{k-1,n_{k-1}}\},
\end{equation*}and 
\begin{equation*}
\mathcal{B}^{p,r}_{k} = \mathrm{span}\{B^{p,r}_{k,1},B^{p,r}_{k,2},B^{p,r}_{k,3},...,B^{p,r}_{k,n_k}\}.
\end{equation*}
The following result  expresses coarse basis functions as the linear combination of fine basis functions.
\begin{proposition}
\label{le:BcBf}
Each coarse basis function $B^{p,r}_{k-1,i}, i = 1,2,...,n_{k-1},$ can be represented as the linear combination of the fine basis functions $\{B^{p,r}_{k,i}, i = 1,2,...,n_k\}$ by the following relation  
\begin{equation}
\label{eq:BcBf}
 \mathcal{B}^{p,r}_{k-1} = G^{p,r}_{k} \mathcal{B}^{p,r}_{k}, \quad \mathrm { i.e., }\quad B^{p,r}_{k-1,i} = \displaystyle \sum_{j=1}^{n_k} g_{ij}B^{p,r}_{k,j}, 
\end{equation}
where $G^{p,r}_{k}=(g_{ij})_{n_{k-1} \times n_k}$, is called the restriction operator from a given fine level to the next coarse level for B-spline basis functions.
\end{proposition}

In the following we explain the formation of transfer operator $G^{p,r}_{k}$ at different levels of mesh and with increasing polynomial degree with the $C^{p-1}$ and $C^0$-continuity. 
\subsubsection{$C^{p-1}$-continuity}

The B-spline basis functions $B^{2,p-1}_{1,i},i=1,2,3,$ and $B^{2,p-1}_{2,i},i=1,2,3,4,$ of degree $p=2$ on knots $E_1 = \{0, 0, 0, 1, 1, 1\}$ and $E_2 = \{0, 0, 0, \frac{1}{2},1, 1, 1\}$, respectively, are defined in section \ref{subsec:ExBp-1}.
Clearly, the total number of coarse and fine basis functions are three ($n_{k-1}=3$) and four ($n_k=4$), respectively. The matrix $G^{2,p-1}_{2}=(g_{ij})_{3 \times 4}$ is given by the following representation of coarse basis functions as the linear combination of fine basis functions.
\begin{empheq}[box=\fbox]{align} \nonumber
B^{2,{p-1}}_{1,1} & = g_{11} B^{2,{p-1}}_{2,1}+  g_{12} B^{2,{p-1}}_{2,2}+  g_{13} B^{2,{p-1}}_{2,3}+  g_{14} B^{2,{p-1}}_{2,4}, \\ \nonumber
B^{2,{p-1}}_{1,2} & = g_{21} B^{2,{p-1}}_{2,1}+  g_{22} B^{2,{p-1}}_{2,2}+  g_{23} B^{2,{p-1}}_{2,3}+  g_{24} B^{2,{p-1}}_{2,4}, \\ \nonumber
B^{2,{p-1}}_{1,3} & = g_{31} B^{2,{p-1}}_{2,1}+  g_{32} B^{2,{p-1}}_{2,2}+  g_{33} B^{2,{p-1}}_{2,3}+  g_{34} B^{2,{p-1}}_{2,4}.
\end{empheq}
Equivalently, it can be written as
\begin{equation*}
  \mathcal{B}^{2,p-1}_{1} = G^{2,p-1}_{2}\mathcal{B}^{2,p-1}_{2},
\end{equation*}
where
 \[  \mathcal{B}^{2,p-1}_{1} =  \left[ \begin{array}{c}
B^{2,{p-1}}_{1,1} \\
B^{2, {p-1}}_{1,2} \\
B^{2, {p-1}}_{1,3} \end{array} \right],
G^{2,p-1}_{2} =  \left[ \begin{array}{cccc}
 g_{11}  &  g_{12}  &  g_{13} & g_{14}  \\
 g_{21}  &  g_{22}  &  g_{23} & g_{24}  \\
 g_{31}  &  g_{32}  &  g_{33} & g_{34}  \end{array} \right],
 \mathcal{B}^{2,p-1}_{2} =  \left[ \begin{array}{c}
B^{2, {p-1}}_{2,1} \\
B^{2, {p-1}}_{2,2} \\
B^{2, {p-1}}_{2,3} \\
B^{2, {p-1}}_{2,4} \end{array} \right].\] 
For the above set of basis functions, $G^{2,p-1}_{2}$ is given by 
\begin{subequations}
\begin{equation}  
G^{2,p-1}_{2} = \frac{1}{4} \left[ \begin{array}{cccc}
4  &  2  & 0 & 0  \\
 0  &   2  & 2  & 0 \\
 0 &  0  & 2  & 4 \end{array} \right].
\end{equation}
Similarly, the coarse basis functions for $B^{2,{p-1}}_{2,i},i=1,2,3,4,$ at level $2$, can be obtained in terms of $B^{2,{p-1}}_{3,i},i=1,2,...,6,$ by the following matrix
\begin{equation}  
G^{2,p-1}_{3} = \frac{1}{4}  \left[ \begin{array}{cccccc}
     4&     2&     0&     0&     0&     0\\
     0&     2&     3&     1&     0&     0\\
     0&     0&     1&     3&     2&     0\\
     0&     0&     0&     0&     2&     4
\end{array} \right].
\end{equation}
In a multilevel setting, the representation of each basis function $B^{2,{p-1}}_{k-1,i}$ at level $k-1$ as the linear combination of the basis functions $B^{2,{p-1}}_{{k},i}$ at level ${k}$ is given by the the following matrix $G^{2,p-1}_{{k}}$, where $k\ge4$.
\begin{equation}  
G^{2,p-1}_{{k}} = \frac{1}{4} \left[ \begin{array}{cccccccccccccccccc}
     4  &   2&     &     &     &      &     &    &     &     &     &     &     &     &     &      &     &      \\
      &    2&     3&     1 &    &     &      &    &     &     &     &     &     &     &     &     &      &     \\
      &    &     1&     3 &    3  &   1&      &    &     &     &     &     &     &    &      &    &      &     \\
         & &     &     &     1  &   3&     3&     1&     &     &     &     &     &     &     &     &      &     \\
     &     &     &     &      &    &     ..&     ..&     ..&     ..&     &     &     &     &     &     &      &    \\
     &     &     &     &      &    &     &     &     ..&     ..&   ..  & ..    &     &     &     &     &      &     \\
         & &     &      &     &    &      &    &     &     &     1&     3&     3&     1&     &      &     &       \\
     &     &     &     &     &     &     &     &     &     &     &     &     1&     3&     3&     1 &     &        \\
       &   &     &      &    &     &      &    &     &     &     &     &     &     &     1&     3 &    2&            \\
        &  &     &       &   &     &     &     &     &     &     &     &    &     &      &      &   2 &    4
 \end{array} \right].
\end{equation}
 The size of the matrix $G^{2,p-1}_{{k}}$ is $(n_{{k-1}}+2) \times (n_{{k}} +2)$, where $n_{{k-1}}$ and $n_{{k}}$ are the number of total knot spans at level $k-1$ and ${k}$, respectively.
\end{subequations}

For higher degree polynomials, the transfer operators can be defined in a similar way. For $p=3$, at level $k=1$ the basis functions $B^{3,p-1}_{1,i},i=1,2,3,4,$ with $C^{p-1}$-continuity  can be represented by the following restriction operator at level $l=2$.
\begin{subequations}
\begin{equation}  
G^{3,p-1}_{2} =\frac{1}{2}  \left[ \begin{array}{ccccc}
     2       &       1   &        0      &        0    &          0\\       
     0       &       1     &       1    &        0    &          0       \\
     0       &       0        &      1     &       1   &         0       \\
     0       &       0        &      0        &      1     &       2       
 \end{array} \right].
\end{equation}
The transfer operator for level $3$ can be written as 
\begin{equation}  
G^{3,p-1}_{3} = \frac{1}{16}\left[ \begin{array}{ccccccc}
        16         &     8           &0&              0              &0&              0              &0       \\
       0          &  8           & 12  &          3            &0      &        0           &   0       \\
       0        &      0       &      4     &       10       &    4        &    0         &     0       \\
       0      &        0    &          0       &       3   &        12        &  8     &      0       \\
       0    &          0 &             0          &    0&              0          & 8&            16       
 \end{array} \right].
\end{equation} 
For all levels ${k}$, where $k\ge4$, we have  
\begin{equation}  
G^{3,p-1}_{{k}} = \frac{1}{16}   \left[ \begin{array}{cccccccccccccccc}

      16     &         8           &&&&&&&&&&&&&&     \\ 
    &           8      &       12      &        3          &&&&&&&&&&&&\\
 &&     4     &       11     &         8    &          2      &&&&&&&&&&    \\
   &&&          2      &        8       &      12      &        8      &        2         &&&&&&&&  \\
&&&&&..& ..&     ..       &     ..     &    ..         &&&&&&   \\
&&&&&& &    ..      &     ..    &    ..          &..&..&&&&   \\
   &&&&&&&&        2     &         8        &     12    &          8     &        2        &&&\\
     &&&&&&&&&&             2     &         8     &         11     &       4         &&    \\
    &&&&&&&&&&&&          3       &       12      &       8           &       \\
    &&&&&&&&&&&&&&           8    &         16      
 \end{array} \right].
\end{equation}
\end{subequations}
The size of the matrix $G^{3,p-1}_{{k}}$ is $(n_{{k-1}}+3) \times (n_{{k}}+3)$.

Finally, we give the transfer operators for $p=4$ with $C^{p-1}$-continuity. For levels $2$ and $3$ the transfer operators are defined as follows:
\begin{subequations}
\begin{equation} 
G^{4,p-1}_{2} = \frac{1}{2} \left[ \begin{array}{cccccc}
       2&              1             & 0 &             0             & 0 &             0        \\
       0   &           1          &    1    &          0          &    0    &          0       \\
       0      &        0       &       1       &       1       &       0       &       0     \\  
       0         &     0    &          0          &    1    &          1          &    0   \\    
       0            &  0 &             0             & 0 &             1             & 2   
 \end{array} \right],
\end{equation}
\begin{equation}  
G^{4,p-1}_{3} = \frac{1}{48} \left[ \begin{array}{cccccccc}
      48             &24&              0             & 0&              0             &0&              0              &0       \\
       0            & 24    &         36          &    9   &           0          &    0  &            0           &   0       \\
       0         &     0        &     12        &     30      &        9       &       0     &         0        &      0       \\
       0      &        0           &   0      &        9           &  30     &        12        &      0     &         0       \\
       0   &           0              &0   &           0              &9   &          36            & 24  &            0       \\
       0&              0              &0&              0              &0&              0             &24&             48       
 \end{array} \right],
\end{equation}
respectively. For levels $k$,  where $k\ge4$, the transfer operator is given by the following
\begin{equation} 
G^{4,p-1}_{{k}} =\frac{1}{48} \left[ \begin{array}{cccccccccccccccc}
      48          &     24           &&&&&&&&&&&&&& \\ 
      &         24     &        36        &      9          &&&&&&&&&&&& \\
    &&          12         &    33   &          20    &          4          &&&&&&&&&& \\
 &&&      6      &       25     &        29      &       15        &      3       &&&&&&&&\\
   &&&&         3    &         15      &       30    &         30       &      15       &       3              &&&&&&\\
       &&&&&&      ..  &        ..    &      ..       &     ..        &&&&&&\\
       &&&&&&      ..  &        ..    &      ..       &     ..        &&&&&&\\
   &&&&&&        3         &    15      &       30        &     30        &     15        &     3              &&&&  \\
      &&&&&&&&             3       &      15     &        29            &  25     &         6       &&&  \\
      &&&&&&&&&&               4          &     20        &     33       &      12          &&  \\
      &&&&&&&&&&&&         9         &    36       &      24             &   \\
    &&&&&&&&&&&&&&         24     &        48  
\end{array} \right],
\end{equation}
\end{subequations}
where the size of the matrix is $(n_{{k-1}}+4) \times (n_{{k}}+4)$. 
\begin{remark}
Since in the  span of an internal basis function of degree $p$ at coarse level, there are $p+2$ full basis functions in the same span at fine level, therefore, any row of $G^{p,p-1}_{{k}}$
can have at most $p+2$ nonzero entries.
\end{remark}

\subsubsection{$C^0$-continuity}

In section \ref{subsec:ExB0}, we explained the explicit representation of $C^0$ continuous B-spline basis functions. The corresponding transfer operators are given in this section. 
The transfer operator $G^{2,0}_{2}$ for $p=2$ with $C^0$-continuity at level $2$ is given by  
\begin{subequations}
\begin{equation}  
G^{2,0}_{2} =\frac{1}{4}  \left[ \begin{array}{ccccc}
     4&     2&     1&     0&     0\\
     0&     2&     2&     2&     0\\
     0&     0&     1&     2&     4
\end{array} \right].
\end{equation}
The operator $G^{2,0}_{{k}}$,  where $k\ge3$ is given by
\begin{equation} 
G^{2,0}_{{k}} =\frac{1}{4} \left[ \begin{array}{ccccccccccccccccccc}
\cline{1-5}
      \multicolumn{1}{|c}     {4}  &  2   &  1   &0&0&  \multicolumn{1}{|c}{}  &&&&&&&&&&&&&\\
     \multicolumn{1}{|c}    { 0}    &  2  &   2  &   2  &0 &  \multicolumn{1}{|c}{}  &&&&&&&&&&&&&\\
\cline{5-9} 
     \multicolumn{1}{|c}    { 0}    &0&   1 &    2  &    \multicolumn{1}{|c|}     {4}   &   2  &   1    &0&   \multicolumn{1}{c|}{0}&  &&&&&&&&&\\
\cline{1-5}
     &&&&  \multicolumn{1}{|c}    { 0}    &   2  &   2   &  2   &     \multicolumn{1}{c|}{0}   &&&&&&&&&&\\
&&&&  \multicolumn{1}{|c}    { 0}    &   0  &   1  &  2   &     \multicolumn{1}{c|}{4}   &&&&&&&&&&\\
\cline{5-9}
     &&&&& &  &  & .. &..&..&&&&&&&&\\
     &&&&&   &      &     &..&..&..&&&&&&&&\\
\cline{11-15}
    &&&&&&&&   && \multicolumn{1}{|c}{4}&  2  &   1  &   0     & \multicolumn{1}{c|}{0}&&&&\\
     &&&&&&&&   && \multicolumn{1}{|c}{0}&  2  &   2  &   2     & \multicolumn{1}{c|}{0}&&&&\\
\cline{15-19}
     &&&&&&&&&&\multicolumn{1}{|c}{0}&0&1& 2&   \multicolumn{1}{|c|}    { 4} &2&1  &   0  &         \multicolumn{1}{c|}{0}\\
\cline{11-15}
      &&&&&&&&&&&&& &   \multicolumn{1}{|c}    { 0} &2&2  &   2  &         \multicolumn{1}{c|}{0}\\
     &&&&&&&&&&&&& &   \multicolumn{1}{|c}    { 0} &0&1  &   2  &         \multicolumn{1}{c|}{4}\\
\cline{15-19}
 \end{array} \right],
\end{equation}
\end{subequations}
 with size $(2n_{{k-1}}+1) \times (2n_{{k}}+1)$. The  matrix $G^{2,0}_{{k}}, k \ge 3,$ has block structure with blocks $G^{2,0}_{{2}}$. The blocks are connected in such a way that if a block ends at $i$th row and $j$th column of $G^{2,0}_{{k}}$ then the next block will start at $(i,j)$th position of $G^{2,0}_{{k}}$ with an overlap of last entry and first entry of the corresponding blocks (which are same).

The transfer operators for $p=3$ with $C^0$-continuity for level $2$ is given by
\begin{equation} 
G^{3,0}_{{2}} = \frac{1}{8}  \left[ \begin{array}{ccccccc}
       8&              4          &    2 &             1           &   0 &             0           &   0 \\      
       0   &           4       &       4    &          3        &      2    &          0        &      0   \\    
       0      &        0    &          2       &       3     &         4       &       4     &         0     \\  
       0         &     0 &              0         &     1 &             2         &     4 &             8       
 \end{array} \right].
\end{equation}
Following the same block structure as in $G^{2,0}_{{k}}$, we can generate $G^{3,0}_{{k}}$,  where $k\ge3$ with size $(3n_{{k-1}}+1) \times (3n_{{k}}+1)$. Finally for $p=4$, we have the following transfer operator for level $2$
\begin{equation} 
G^{4,0}_{2} =\frac{1}{16}  \left[ \begin{array}{ccccccccc}
     16&              8             & 4 &             2             & 1&              0             & 0 &             0              &0  \\     
       0   &           8          &    8    &          6          &    4   &           2          &    0    &          0           &   0    \\   
       0      &        0       &       4       &       6       &       6      &        6       &       4       &       0       &       0      \\ 
       0         &     0    &          0          &    2    &          4         &     6    &          8          &    8    &          0       \\
       0            &  0 &             0             & 0 &             1            &  2 &             4             & 8 &            16       
 \end{array} \right],
\end{equation}
Similarly, repeating these blocks as in previous cases, we can generate $G^{4,0}_{{k}}$,  where $k\ge3$  with size $(4n_{{k-1}}+1) \times (4n_{{k}}+1)$.

\begin{remark}
 Note that the transfer operators are defined for one dimensional B-splines. For two- and three-dimensions, we take tensor product of these operators.
\end{remark}

\subsection{Multilevel NURBS}

This section presents the procedure for constructing NURBS multilevel spaces in a simplified manner. Since NURBS are generated from B-splines, its natural to construct NURBS transfer operators from B-splines transfer operators. For a two level setting, let $\mathcal{N}^{p,r}_{k-1}$ and $\mathcal{N}^{p,r}_{k}$ denote the NURBS spaces at coarse and fine level, respectively. Let $\{N^{p,r}_{k-1,i}, i = 1,2,...,n_{k-1}\}$ and $\{N^{p,r}_{k,i}, i = 1,2,...,n_k\}$ be the set of basis functions for coarse and fine space, respectively, i.e.
\begin{equation*}
\mathcal{N}^{p,r}_{k-1} = \mathrm{span}\{N^{p,r}_{k-1,1},N^{p,r}_{k-1,2},N^{p,r}_{k-1,3},...,N^{p,r}_{k-1,n_{k-1}}\},
\end{equation*}and 
\begin{equation*}
\mathcal{N}^{p,r}_{k} = \mathrm{span}\{N^{p,r}_{k,1},N^{p,r}_{k,2},N^{p,r}_{k,3},...,N^{p,r}_{k,n_k}\}.
\end{equation*}
Note that, a relation similar to Proposition \ref{le:BcBf} also holds for NURBS basis functions, i.e., we have 
\begin{equation}
\label{eq:NcNf}
 \mathcal{N}^{p,r}_{k-1}= R^{p,r}_{k} \mathcal{N}^{p,r}_{k}, \quad \text{ i.e., } \quad N^{p,r}_{k-1,i} = \displaystyle \sum_{j=1}^{n_k} r_{ij}N^{p,r}_{k,j}, \quad \forall i = 1,2,3,...,n_{k-1},
\end{equation}
where $R^{p,r}_{k}=(r_{ij})_{n_{k-1} \times n_k}$, is  restriction operator with respect to NURBS basis functions. As NURBS are formed from B-splines and weights, $R^{p,r}_{k}$ can be obtained from $G^{p,r}_{k}$ and weights. Using the definition of NURBS and \eqref{eq:NcNf}, we have 
\begin{equation}
\label{eq:nurbscf}
\begin{split}
 \frac{ w_i^{k-1} B^{p,r}_{k-1,i}}{\displaystyle  \sum_{i'=1}^{n_{k-1}} w_{i'}^{k-1} B^{p,r}_{k-1,i'}} &=  \displaystyle \sum_{j=1}^{n_k} r_{ij} \frac{ w_j^{k} B^{p,r}_{k,j}}{\displaystyle  \sum_{j'=1}^{n_k} w_{j'}^k B^{p,r}_{k,j'}},\quad \forall i = 1,2,3,...,n_{k-1},
\end{split}
\end{equation}
where $w_i^{k-1}, i=1,2,3,...,n_{k-1},$ and $w_j^k, j = 1,2,3,...,n_k,$ are the weights for coarse space and fine space, respectively. Note that the weight function  $\displaystyle  \sum_{i=1}^{n} w_i B_i$ does not change its value with respect to refinements, i.e., we have
\begin{equation}
\label{eq:weightref}
  \displaystyle  \sum_{i=1}^{n_{k-1}} w_i^{k-1} B^{p,r}_{k-1,i} = \displaystyle  \sum_{j=1}^{n_k} w_j^k B^{p,r}_{k,j},
\end{equation}
which is an important result from the refinement point of view.  Now using \eqref{eq:weightref}, from \eqref{eq:nurbscf} we get
\begin{equation*}
 \hspace{1cm}\displaystyle { w_i^{k-1} B^{p,r}_{k-1,i}} =  \displaystyle \sum_{j=1}^{n_k} r_{ij} { w_j^k B^{p,r}_{k,j}},
\end{equation*}
and thus
\begin{equation}
\label{eq:BcRBf}
\displaystyle {B^{p,r}_{k-1,i}}  =  \displaystyle \sum_{j=1}^{n_k}\frac{ r_{ij} w_j^k }{ w_i^{k-1} } {B^{p,r}_{k,j}}.
\end{equation}
Comparing the coefficients of $B^{p,r}_{k,j}$ in  \eqref{eq:BcBf} and \eqref{eq:BcRBf}, we get 
\begin{equation}
\label{eq:BcRRBf}
 \frac{r_{ij} w_j^k }{ w_i^{k-1} }=g_{ij} \Longrightarrow  r_{ij}= \frac{ w_i^{k-1} g_{ij}}{w_j^k }.
\end{equation}
This can be equivalently written as 
\begin{equation}
\label{eq:RrelationG}
R^{p,r}_k={W_I^{k-1}G^{p,r}_k}{\left({W_I^k}\right)^{-1}},
\end{equation}
where $W^{k-1}_I$ and $W^k_I$ are the diagonal matrices corresponding to the weights at the coarse level and the fine level, respectively, and defined as follows

\[
W_I^{k-1}= \left[ \begin{array}{ccccccccc}
w^{k-1}_1 & &  & & & & & \\
 & w^{k-1}_2& &  & & & & \\
 & & & ..& &  & & \\
 & & & &..&  & & \\
 & & & & &&w^{k-1}_{n_{k-1}-1} & \\
 & & & & & & &w^{k-1}_{n_{k-1}} 
 \end{array} \right],\]
\[
W_I^k = \left[ \begin{array}{ccccccccc}
w^k_{1}  & &  & & & & & \\
 & w^k_{2} & &  & & & & \\
 & & & ..& &  & & \\
 & & & &..&  & & \\
 & & & & &&w^k_{n_k-1}  & \\
 & & & & & & &w^k_{n_k} 
 \end{array} \right].\] 
The equation \eqref{eq:RrelationG} gives us the NURBS operators using B-splines transfer operators and weights at coarse and fine levels. From \eqref{eq:weightref} we can also obtain the procedure to refine the weights as follows. We have 
\begin{equation*}
\displaystyle  \sum_{i=1}^{n_{k-1}} w_i^{k-1} B^{p,r}_{k-1,i}  = \displaystyle  \sum_{j=1}^{n_k} w_j^k B^{p,r}_{k,j},
\end{equation*}
which implies
\begin{equation*}
\displaystyle  \sum_{i=1}^{n_{k-1}} w_i^{k-1} \displaystyle \sum_{j=1}^{n_k} g_{ij}B^{p,r}_{k,j}  = \displaystyle  \sum_{j=1}^{n_k} w_j^k B^{p,r}_{k,j}. 
\end{equation*}
Comparing the coefficients of  $B^{p,r}_{k,j}$ from both the sides, we get
\begin{equation}
 w_j^k = \displaystyle \sum_{i=1}^{n_{k-1}} w_i^{k-1} g_{ij}\quad \text{ for } j= 1,2,...,n_k.
\end{equation}
Equivalently, this can be written in matrix form as follows
\begin{equation}
 W^k = \left({G^{p,r}_k}\right)^TW^{k-1},
\end{equation}
where 
\[W^k=
\left[ \begin{array}{c}
w_1^k\\
w_2^k\\
:\\
:\\
w_{n_k-1}^k\\
w_{n_k}^k
 \end{array} \right], W^{k-1} =
\left[ \begin{array}{c}
w_1^{k-1}\\
w_2^{k-1}\\
:\\
:\\
w_{n_{k-1}-1}^{k-1}\\
w_{n_{k-1}}^{k-1}
 \end{array} \right].\] 
Using above, now we can write the NURBS operators in terms of B-spline operator and weights only at coarse level. From \eqref{eq:BcRRBf}, we get
\begin{equation}
 r_{ij}  = \displaystyle \frac{{g_{ij} w_i^{k-1}}}{{\displaystyle \sum_{i=1}^{n_{k-1}} w_i^{k-1} g_{ij}}}.
\end{equation}
In matrix form this can be written as 
\begin{equation}
R^{p,r}_k={W^{k-1}_I G^{p,r}_k }{\left(\mathrm{diag}\left(\left({G^{p,r}_k}\right)^TW^{k-1}\right)\right)}^{-1}.
\end{equation}
\begin{remark}
 The operators $G^{p,r}_k$ and $R^{p,r}_k$ can also be used in constructing restriction operators in multigrid methods, see e.g., \cite{GahalautKT-12}. 
\end{remark}
\begin{remark}
In practice, these operators are constructed once for all levels and stored in sparse matrix format.
\end{remark}

\section{AMLI Methods}
\label{sec:amli_methods}

In this section we present the basic principle of AMLI methods.
In what follows we will denote by $M^{(k)}$ a preconditioner for the stiffness matrix $A^{(k)}$ corresponding to level $k$. We will also make use of the corresponding hierarchical matrix $\hat{A}^{(k)}$, which is related to $A^{(k)}$ via a two-level hierarchical basis (HB) transformation $J^{(k)}$, i.e.,
\begin{equation}\label{NBHB}
\hat{A}^{(k)} = J^{(k)} A^{(k)} (J^{(k)})^T .
\end{equation}
The transformation matrix $J^{(k)}$ specifies the space splitting, which will be described in detail in Section \ref{sec:hier_space}. By  $A^{(k)}_{ij}$ and $\hat{A}^{(k)}_{ij}$, $1 \le i,j \le 2$, we denote the blocks of $A^{(k)}$ and $\hat{A}^{(k)}$ that correspond to the fine-coarse partitioning of degrees of freedom where the degrees of freedom associated with the coarse mesh are numbered last.

\medskip

The aim is to build a multilevel preconditioner $M^{(L)}$ for the coefficient matrix $A^{(L)}$ at the finest level that has a uniformly bounded (relative) condition number
\[
\varkappa({M^{(L)}}^{-1} A^{(L)}) = {\mathcal O}(1),
\]
and an optimal computational complexity, that is, linear in the number of degrees of freedom $N_{L}$ at the finest level. In order to achieve this goal hierarchical basis methods can be combined with various types of stabilization techniques.

\medskip

One particular purely algebraic stabilization technique is the so-called algebraic  multilevel iteration (AMLI) method, where a specially constructed matrix polynomial $p^{(k)}$ of degree $\nu_{k}$ can be employed at some (or all) levels $k$. 

\medskip

We have the following two-level hierarchical basis representation at level $k$
\begin{equation}
\label{eq:2x2}
\hat{A}^{(k)} =
\begin{bmatrix}
\hat{A}_{11}^{(k)}&\hat{A}_{12}^{(k)}\\
\hat{A}_{21}^{(k)}&\hat{A}_{22}^{(k)}\end{bmatrix},
\end{equation}
where $\hat{A}_{22}^{(k)}=A^{({k-1})}$ is the coarse-level stiffness matrix.
Starting at level $l=1$ (associated with the coarsest mesh), on which a complete LU factorization of the matrix $A^{(1)}$ is performed, we define
\begin{equation}\label{ml}
M^{(1)}:=A^{(1)}.
\end{equation}
Given the preconditioner $M^{({k-1})}$ at level ${k-1}$, the preconditioner $M^{(k)}$ at level $k$ is then defined by
\begin{equation}\label{mk}
M^{(k)}:=L^{(k)} U^{(k)},
\end{equation}
where
\begin{equation}\label{rklk}
L^{(k)}:=\left[ \begin{array}{cc}
C_{11}^{(k)} & 0 \\
\hat{A}_{21}^{(k)} & C_{22}^{(k)}
\end{array} \right], \quad
U^{(k)}:=\left[ \begin{array}{cc}
I & {C_{11}^{(k)}}^{-1} \hat{A}_{12}^{(k)} \\
0 & I
\end{array} \right] .
\end{equation}
Here $C_{11}^{(k)}$ is a preconditioner for the pivot block $A_{11}^{(k)}$, and
\begin{equation}\label{zk}
C_{22}^{(k)} := {A}^{({k-1})} \left( I - p^{(k)} ( {M^{({k-1})}}^{-1} A^{({k-1})} ) \right)^{-1}
\end{equation}
\begin{equation}\label{pk}
0 \le p^{(k)} (t) \le 1, \quad 0 \le t \le 1, \quad p^{(k)} (0) = 1 .
\end{equation}
It is easily seen that (\ref{zk}) is equivalent to
\begin{equation}\label{zk1}
{C_{22}^{(k)}}^{-1} =
{M^{({k-1})}}^{-1} q^{(k)} ( A^{({k-1})} {M^{({k-1})}}^{-1} ),
\end{equation}
where the polynomial $q^{(k)}$ is given by
\begin{equation}\label{q}
q^{(k)}(x)=\frac{1-p^{(k)}(x)}{x}.
\end{equation}
We note that the multilevel preconditioner defined via (\ref{mk}) is getting close to a two-level method when $q^{(k)}(x)$ closely approximates $1/x$, in which case ${C_{22}^{(k)}}^{-1} \approx {A^{({k-1})}}^{-1}$. In order to construct an efficient multilevel method, the action of ${C_{22}^{(k)}}^{-1}$ on an arbitrary vector should be much cheaper to compute (in terms of the number of arithmetic operations) than the action of ${A^{({k-1})}}^{-1}$. Optimal order solution algorithms typically require that the arithmetic work for one application of ${C_{22}^{(k)}}^{-1}$ is of the order ${\mathcal O}(N_{{k-1}})$, where $N_{{k-1}}$ denotes the number of unknowns at level ${k-1}$. 

\medskip

It is well known from the theory introduced in \cite{AxelssonV-89,AxelssonV-90} that a properly shifted and scaled Chebyshev polynomial $p^{(k)} := p_{\nu_{k}}$ of degree $\nu_{k}$ can be used to stabilize the condition number of ${M^{(k)}}^{-1} A^{(k)}$ (and thus obtain optimal order computational complexity). Other polynomials such as the best polynomial approximation of $1/x$ in uniform norm also qualify for stabilization, see, e.g., \cite{KrausVZ-12}.
Alternatively, in the nonlinear AMLI method, see, e.g., \cite{AxelssonV-94}, a few inner flexible conjugate gradient (FCG) type iterations (for the FCG algorithm, see also \cite{Notay-00}) are performed in order to improve (or freeze) the residual reduction factor of the outer FCG iteration. In general, the resulting nonlinear (variable step) multilevel preconditioning method is of comparable efficiency, and, because its realization does not rely on any spectral bounds, is easier to implement than the linear AMLI method (based on a stabilization polynomial). For a convergence analysis of nonlinear AMLI see, e.g., \cite{Kraus-02, Vassilevski-08}.
%

Typically, the iterative solution process is of optimal order of computational complexity if the degree $\nu_{k} =\nu$ of the matrix polynomial (or alternatively, the number of inner iterations for nonlinear AMLI) at level $k$ satisfies the optimality condition
\begin{alignat}{1}
1/\sqrt{(1- \gamma^{2})} < \, \nu  < \,\tau,
\label{eq:CondOpt}
\end{alignat}
where $ \tau \approx  \tau_{k} = {N_{k}} / {N_{{k-1}}}$ denotes the reduction factor of the number of degrees of freedom, and $\gamma$ denotes the constant in the strengthened Cauchy-Bunyakowski-Schwarz (CBS) inequality. In case of (standard) full coarsening the value of $\tau$ is approximately $4$ and $8$ for two- and three-dimensional problems, respectively. For a more detailed discussion of AMLI methods, including implementation issues, see, e.g., \cite{KrausBook, Vassilevski-08}.

\medskip

\begin{remark}
The AMLI algorithm has originally been introduced and studied in the multiplicative form
\eqref{mk}--\eqref{rklk}, see \cite{AxelssonV-89,AxelssonV-90}.
However, it is also possible to construct the preconditioner in additive form, which
is defined as follows
\begin{equation}
\label{mk_a}
M^{(k)}_{A} := \left[ \begin{array}{cc}
C_{11}^{(k)} & 0 \\
0 & C_{22}^{(k)}
\end{array} \right].
\end{equation}
In this case
the optimality condition for the polynomial degree (or the number of inner iterations
at level $k-1$ induced by one nonlinear AMLI cycle at level $k$) reads
\begin{alignat}{1}
\sqrt{(1 + \gamma)/(1- \gamma)} < \, \nu  < \,\tau.
\label{eq:CondOpt_Add}
\end{alignat}
For details, see \cite{Axelsson-99}.
\end{remark}
%

\section{Hierarchical Spaces}
\label{sec:hier_space}

\subsection{Construction}\label{sec:hier_space_constr}

Hierarchical basis techniques, in the present context, serve the purpose to decompose 
the finite-dimensional spaces of the B-spline (NURBS) basis functions into a coarse space
and its hierarchical complement. In the AMLI framework it is crucial that the angle between
theses two subspaces is uniformly bounded with respect to the mesh size. This issue will be
addressed in Section~\ref{sec:num_cbs}.

From Section~\ref{sec:amli_methods} we recall the following two-level hierarchical basis
representation for stiffness matrix at fine level
\begin{equation}
\label{eq:2x2_again}
\hat{A}^{(k)} =
\begin{bmatrix}
\hat{A}_{11}^{(k)}&\hat{A}_{12}^{(k)}\\
\hat{A}_{21}^{(k)}&\hat{A}_{22}^{(k)}
\end{bmatrix}
= \begin{bmatrix}
\hat{A}_{11}^{(k)}&\hat{A}_{12}^{(k)}\\
\hat{A}_{21}^{(k)}&A^{(k-1)}
\end{bmatrix},
\end{equation} where $\hat{A}_{22}^{(k)}$ represents the matrix corresponding to coarse basis functions and  $\hat{A}_{11}^{(k)}$ represents the matrix corresponding to its hierarchical complement, and $1 \le k \le L$. From Section \ref{sec:MulSplines} recall that, for B-splines we have the following transformations
\begin{equation}
 \hat{A}_{22}^{(k)} =G^{p,r}_k A^k  {(G^{p,r}_k)}^T,
\end{equation} respectively.
For hierarchical complementary spaces, let $T^{p,r}_k$ be the matrix  such that
\begin{equation}
\hat{A}_{11}^{(k)} = T^{p,r}_k A^k  {(T^{p,r}_k)}^T.
\end{equation}
Here the matrix $T^{p,r}_k$ is a hierarchical complementary transfer operator, which transfers fine basis functions to a set of hierarchical complementary basis functions.
The remaining two blocks of the hierarchical matrix $\hat{A}^{(k)}$ can be obtained by the following relations
\begin{equation}
 \hat{A}_{12}^{(k)} =T^{p,r}_k A^{(k)}  {(G^{p,r}_k)}^T, \quad 
 \hat{A}_{21}^{(k)} =G^{p,r}_k A^{(k)}  {(T^{p,r}_k)}^T.
\end{equation} 
Hence, the transformation matrix $J^{(k)}$ in~\eqref{NBHB} has the form
$$
J^{(k)} =
\begin{bmatrix}
T^{p,r}_k \\
G^{p,r}_k
\end{bmatrix}^T.
$$
Note that similar results hold for $R^{p,r}_k $.

\medskip

To construct  $T^{p,r}_k$ efficiently, the following points are important.
\begin{enumerate}
\item The basis for hierarchical complementary space should be locally supported. In other words, the block $\hat{A}_{11}^{(k)}$ should be sparse. 
\item  The condition number of $\hat{A}_{11}^{(k)}$ should be independent of the mesh size.
\item  The CBS constant $\gamma$, see \eqref{eq:gamma}, should be bounded away from one,
i.e. the minimum generalized eigenvalue of the Schur complement  with respect to  $\hat{A}_{22}^{(k)}$ block should be greater than $1/4$ for $\nu=2$ and $1/9$ for $\nu=3$.
\end{enumerate}
The construction of $T^{p,r}_k$, based on the linear combination of fine basis functions, is not unique. Based on the above mentioned guidelines, a representation of a complementary basis function should not involve several fine basis functions because this, in general, will cause more entries in
$T^{p,r}_k$.
Based on our extensive study with different choices of linear combinations satisfying the above requirements, we present two choices of $T^{p,r}_k$, for $p=2,3,4$ and for
$C^{p-1}$ and $C^0$ continuity.

\subsubsection{$C^{p-1}$-continuity}

For the first choice of $T^{p,r}_k$ we have the following matrix representation of the
hierarchical complementary space for $p=2$ with $C^{p-1}$-continuity.
\[T^{2,p-1}_{{k}}= \left[ \hspace{-3mm}\begin{array}{cccccccccccccccccccc}
\cline{2-7}
&\multicolumn{1}{|c}{ 0} &     1       &      -1        &      0      &0&0  & \multicolumn{1}{|c} {} &&&&&&&&&&&&\\
& \multicolumn{1}{|c}{0}  &   0 &            0       &       1      &       -1       &  {   0}         & \multicolumn{1}{|c}{}&&&&&&&&&&&&\\
\cline{2-11}
&   &   &&&        \multicolumn{1}{|c}{  0}    &        1         &    -1     &         0           &0&0& \multicolumn{1}{|c} {}&&&&&&&&\\
 &&&&&    \multicolumn{1}{|c}{ 0} &0&        0     &       1       &      -1       &     0           & \multicolumn{1}{|c} {}&&&&&&&&\\
\cline{6-11}
& &&&&&&&&    &          &     ..   &          ..              &&&&&&\\
& &&&&&&&&   &         &     ..   &          ..              &&&&&&\\
\cline{14-19}
  &    &&&&&&&&&&  &&        \multicolumn{1}{|c}{ 0}       &       1     &        -1        &     0       	&0& 0& \multicolumn{1}{|c} {}  \\
  &     &&&&&&&&&&&&  \multicolumn{1}{|c}{0  } &0&          0     &         1      & 	-1 &0& \multicolumn{1}{|c} {}\\
\cline{14-19}
 \end{array} \hspace{-3mm} \right].\] The above matrix has the block structure with blocks, say $M^{2,p-1}_{1}$.
The blocks are connected in such a way that if a block ends at $i$th row and $j$th column of $T^{2,p-1}_{{k}}$ then the next block will start at $(i+1,j-1)$th position of $T^{2,p-1}_{{k}}$. 
In general, for $p=2,3,4,$ we write the following block form of  $T^{p,p-1}_{{k}}$ with blocks $M^{p,p-1}_{1}$
\begin{equation}
\label{eq:TFCp}
T^{p,p-1}_{{k}} = \left[ \begin{array}{ccccc}
M^{p,p-1}_{1} &  &&&\\
&M^{p,p-1}_{1} &  &&\\
&&..&  &\\
&&&M^{p,p-1}_{1} &\\
&&&&M^{p,p-1}_{1}
 \end{array} \right],
\end{equation}
where 
\[M_{1}^{2,p-1} = \left[ \begin{array}{cccccc}
   0  &     1 &  -1    &0&0 &0\\
  0  &0&0 &     1 &  -1    & 0
 \end{array} \right],\] 
\[M_{1}^{3,p-1} = \left[ \begin{array}{ccccccc}
0     &    -1/2      &      3/4      &     -1/2       &    0   &        0          &    0\\
0    &         0       &      0     &      -1/2    &        3/4    &       -1/2     &      0
 \end{array} \right],\]
and
\[M_{1}^{4,p-1} = \left[ \begin{array}{ccccccccc}
0    &         1/2    &       -1     &         1      &       -1/2      &     0   &         0          &    0  \\
0      &          0   &          0    &         1/2  &         -1   &           1     &        -1/2       &     0
 \end{array} \right],\] respectively. 
The blocks are connected as follows; if a block ends at $i$th row and $j$th column
of $T^{p,p-1}_{{k}}$ then the next block will start at $(i+1,j-(p-1))$th position of $T^{p,p-1}_{{k}}$.

For the second choice of $T^{p,r}_k$ we give the following block matrix.
\begin{equation}
\label{eq:TSCp}
T^{p,p-1}_{{k}} = \left[ \begin{array}{ccccc}
M^{p,p-1}_{2} &  &&&\\
&M^{p,p-1}_{2} &  &&\\
&&..&  &\\
&&&M^{p,p-1}_{2} &\\
&&&&M^{p,p-1}_{2}
 \end{array} \right],
\end{equation}
where the blocks $M_{2}^{p,p-1}$ are given by 
\[M_{2}^{2,p-1} = \left[ \begin{array}{cccccc}
   -1/2  &     1 &  -1    & 1/2 &0 &0\\
  0  &0&-1/2  &     1 &  -1    &1/2
 \end{array} \right],\] 
\[M_{2}^{3,p-1} = \left[ \begin{array}{ccccccc}
1/8       &    -1/2      &      3/4      &     -1/2       &     1/8    &        0          &    0\\
0    &         0       &       1/8     &      -1/2    &        3/4    &       -1/2     &       1/8 
 \end{array} \right],\]
and
\[M_{2}^{4,p-1} = \left[ \begin{array}{ccccccccc}
1/4      &         1/2    &       -1     &         1      &       -1/2      &      -1/4   &         0          &    0  \\
0      &          0   &          1/4     &         1/2  &         -1   &           1     &        -1/2       &     -1/4 
\end{array} \right],\] respectively, and the blocks are connected
in the same way as before.

\subsubsection{$C^{0}$-continuity}

For the first choice of $T^{p,r}_k$ with $C^0$-continuity, we give the following matrix representation of the hierarchical complementary spaces. For $p=2$, we have
\[T^{2,0}_{{k}}= \left[ \begin{array}{cccccccccccccc}
\cline{1-5}
\multicolumn{1}{|c}{ 0} &      1 &  -1/4     &0&\multicolumn{1}{c|}{ 0}&&&&&&&&&\\
\multicolumn{1}{|c}{ 0}  &0&      1 &  -1/4    &\multicolumn{1}{c|}{ 0}&&&&&&&&&\\
\cline{1-9}
       &&&&\multicolumn{1}{|c}{ 0}&      1 &  -1/4      &0&\multicolumn{1}{c|}{ 0}&&&&&\\
     &&&&\multicolumn{1}{|c}{ 0}&0&     1 &  -1/4     &\multicolumn{1}{c|}{ 0}&&&&&\\
\cline{5-9}
&&&& & &&      &  ..    & ..&&&&\\
&&&& & &&    &..  & ..&&&&\\
\cline{10-14}
     &&&&&&&&&\multicolumn{1}{|c}{ 0}&      1 &  -1/4    &0&\multicolumn{1}{c|}{ 0}\\
     &&&&&&&&&\multicolumn{1}{|c}{ 0}&0&      1 &  -1/4    &\multicolumn{1}{c|}{ 0}\\
\cline{10-14}
 \end{array} \right].\] 
The above matrix has a block structure and the blocks are connected as follows; if a block ends
at the $i$th row and $j$th column of $T^{2,0}_{{k}}$ then the next block will start at $(i+1,j)$th position of $T^{2,0}_{{k}}$. In general, for $p=2,3,4,$ we can use the following hierarchical complementary operators
\begin{equation}
\label{eq:TFC0}
T^{p,0}_{{k}} = \left[ \begin{array}{ccccc}
M_{1}^{p,0} &  &&&\\
&M_{1}^{p,0} &  &&\\
&&..&  &\\
&&&M_{1}^{p,0} &\\
&&&&M_{1}^{p,0}  
 \end{array} \right],
\end{equation}
where 
\[M_{1}^{2,0} = \left[ \begin{array}{cccccc}
   0  &     1 &  -1/4    &0&0 &0\\
  0  &0&0  &     -1/4 &  1    & 0
 \end{array} \right],\] 
\[M_{1}^{3,0} = \left[ \begin{array}{ccccccc}
       0      &       1     &      -1     &      0           &   0 &            0     &         0       \\
       0      &        0   &         0     &      1/2  &      -1/2     &       0     &         0       \\
       0       &       0      &       0          &    0           &   1         &   -1     &      0       
 \end{array} \right],\] 
and
\[M_{1}^{4,0}= \left[ \begin{array}{ccccccccc}
       0     &        -2/3      &     5/4&            0     &      0 &             0             & 0 &             0           &   0 \\      
       0        &      0      &       -2/3 &           5/4      &      0  &          0          &    0    &          0        &      0   \\    
       0           &   0   &           0       &       0    &         0      &      5/4        &     -2/3     &       0     &         0     \\  
       0              &0&              0          &    0 &             0            &    0  &         5/4           &  -2/3&            0       
 \end{array} \right],\] respectively.

The second choice of $T^{p,r}_k$ for $C^0$ continuous basis functions is obtained by choosing
the following block matrix 
\begin{equation}
\label{eq:TSC0}
T^{p,0}_{{k}} = \left[ \begin{array}{ccccc}
M_{2}^{p,0} &  &&&\\
&M_{2}^{p,0} &  &&\\
&&..&  &\\
&&&M_{2}^{p,0} &\\
&&&&M_{2}^{p,0}  
 \end{array} \right],
\end{equation}
where
\[M_{2}^{2,0}  = \left[ \begin{array}{ccccc}
   -1/4 &     1 &  -1/4    &0&0\\
  0  &0&     -1/4  &  1   & -1/4
\end{array} \right],\] 
\[M_{2}^{3,0}  = \left[ \begin{array}{ccccccc}
       0      &       -1/2      &      1/2      &      0           &   0 &            0     &         0       \\
       0      &        0   &          -1/4      &      1/10    &      -1/4     &       0     &         0       \\
       0       &       0      &       0          &    0           &   1/2        &   -1/2      &      0       
 \end{array} \right],\] 
and
\[M_{2}^{4,0} = \left[ \begin{array}{ccccccccc}
       0     &        -5/9      &      1&             -5/9      &      0 &             0             & 0 &             0           &   0 \\      
       0        &      0      &       -5/9 &           1      &       -5/9  &          0          &    0    &          0        &      0   \\    
       0           &   0   &           0       &       0    &         -5/9      &      1        &     -5/9     &       0     &         0     \\  
       0              &0&              0          &    0 &             0            & -5/9   &         1           &  -5/9&            0       
 \end{array} \right],\]
respectively, and the blocks are connected
in the same way as before.

\begin{remark}
\label{rem:TensorOperator}
All the above operators are defined for one space dimension. The higher dimensional operators are obtained via tensor products.
\end{remark}

\subsection{Quality Assessment}
\label{sec:num_cbs}

The construction of optimal preconditioners in the framework of AMLI methods is
based upon a theory in which the constant $\gamma$ in the strengthened Cauchy-Bunyakowski-Schwarz (CBS) inequality plays a key role. The CBS constant measures the cosine of the abstract
angle between the coarse space and its hierarchical complementary space. The general idea is to construct a
proper splitting by means of a hierarchical basis transformation.

In the hierarchical bases context we denote by $V_1$ and $V_2$ subspaces of the  space $V_h$ .
The space $V_2$ is spanned by the coarse-space basis functions and  $V_1$ is the hierarchical complement of  $V_2$ in  $V_h$,
i.e.,  $V_h$ is a direct sum of  $V_1$ and  $V_2$:
\begin{equation*}
 V_h = V_1 \oplus V_2.
\end{equation*}
Let $v_i \in V_i, i=1,2$. The CBS constant measures the strength of the off-diagonal blocks in relation to the diagonal blocks (see, \eqref{eq:2x2}) and can be defined as the minimal $\gamma$ satisfying the strengthened CBS inequality
\begin{equation}
\label{eq:gamma}
 |v_1^T \hat{A}_{12} v_2| \le \gamma \left\{ (v_1^T \hat{A}_{11}v_1)  (v_2^T \hat{A}_{22} v_2) \right\}^{1/2}.
\end{equation}
A detailed exposition of the role of this constant can be found in \cite{EijkhoutV-91}.

In finite element context, the CBS constant can be estimated locally for various
discretizations and hierarchical transformations. Most of the work on this topic
has been conducted for linear conforming and non-conforming elements, see \cite{KrausBook}
and the references therein, and only very few results exist up to now for quadratic and
higher-order elements, see e.g.~\cite{AxelssonBlaheta-04,EGF-11,KrausLM-12,LymberyMargenov-12}.

The basic idea is as follows. Let us assume that
\begin{equation}
 \hat{A} = \sum_{E \in \mathcal{E}} R_E^T A_E R_E, \quad v = \sum_{E \in \mathcal{E}} R_E^T v_E,
\end{equation}
where $A_E$ are symmetric positive semidefinite local matrices (macro element matrices), $\mathcal{E}$ is some index set, and the summation is understood as assembling. The global splitting naturally induces the two-by-two block representation of the local matrix $A_E$ and the related vector $v_E$, namely,
\[A_E = \left[ \begin{array}{cc}
   A_{E:11}  & A_{E:12}   \\
  A_{E:21}  &  A_{E:22}
 \end{array} \right], \quad v_E = \left[ \begin{array}{c}
   v_{E:1} \\
  v_{E:2} 
 \end{array} \right].\]
Then the local CBS constant $\gamma_E$ corresponding to $A_E$ satisfies the inequality
\begin{equation}
 |v_{E:1}^T A_{E:12} v_{E:2}| \le \gamma_E \left\{ (v_{E:1}^T A_{E:11} v_{E:1})  (v_{E:2}^T A_{E:22} v_{E:2})  \right\}^{1/2}.
\end{equation}
As it is shown in \cite{KrausBook}, the relation between global $\gamma$ and local $\gamma_E$ is given by
\begin{equation}
 \gamma \le \mathrm{max}_{E\in\mathcal{E}} \gamma_E <1.
\end{equation}
\begin{figure}[p]
\caption{B-spline basis functions for $p=2$ on a unit interval with $8$ subdivisions.
The pictures from top to bottom represent basis functions at fine level, coarse level,
for the hierarchical complement of the coarse level, and for the direct sum of the coarse
space and its hierarchical complement respectively.}
\label{fig:HSpace}
 \includegraphics[height=7.5cm,width=15cm]{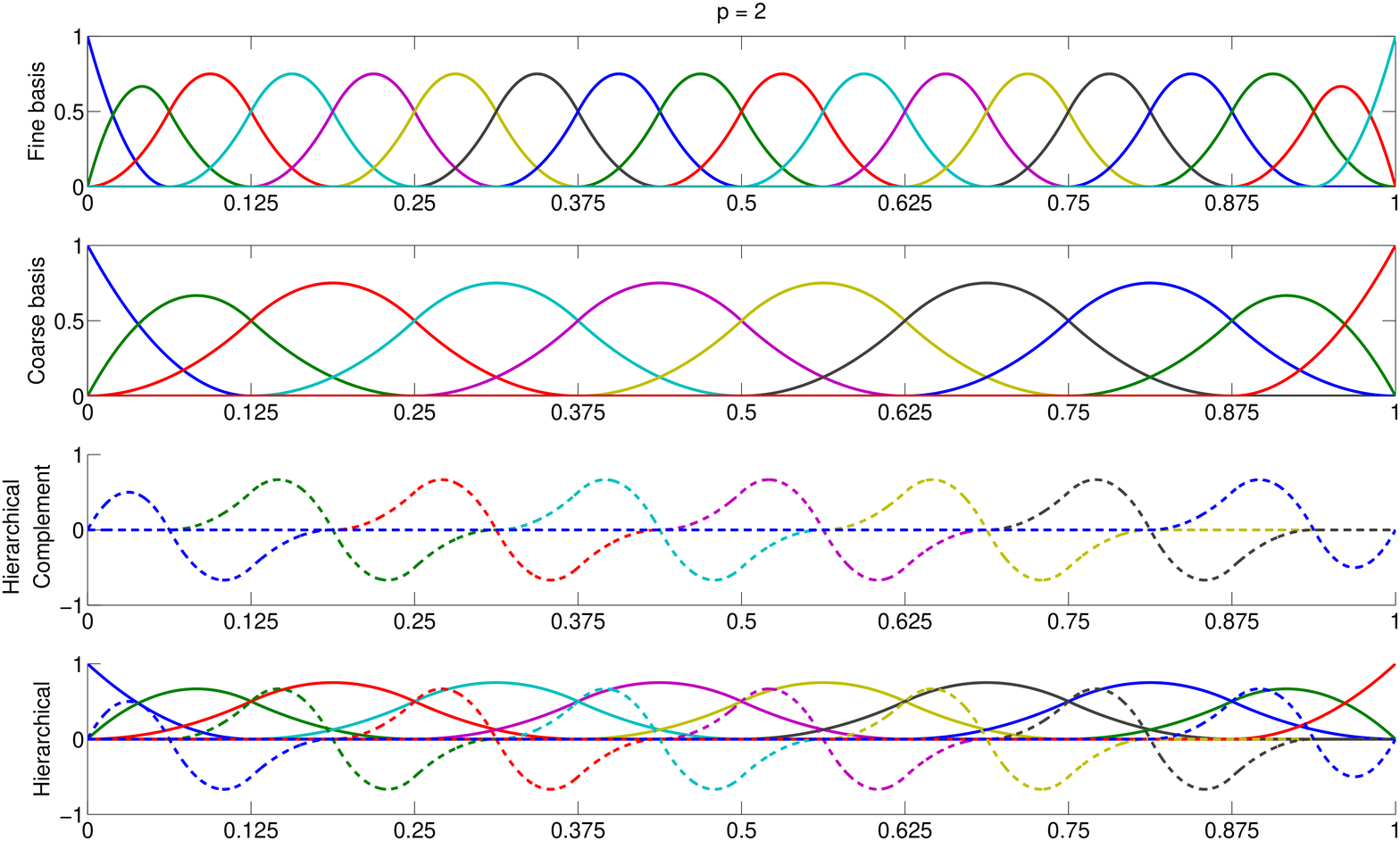}
\end{figure}
\begin{figure}[p]
\caption{Dimension mismatch of basis function for macro element in fine space and its corresponding hierarchical space.  For a given macro element the coarse space and its hierarchical complementary
space have $3$ basis functions each, which results in $6$ basis functions in the hierarchical space, whereas there are only $4$ basis functions for the same macro element in the fine space.\vspace{3mm}}
\label{fig:Dimension}
 \includegraphics[height=7.5cm,width=15cm]{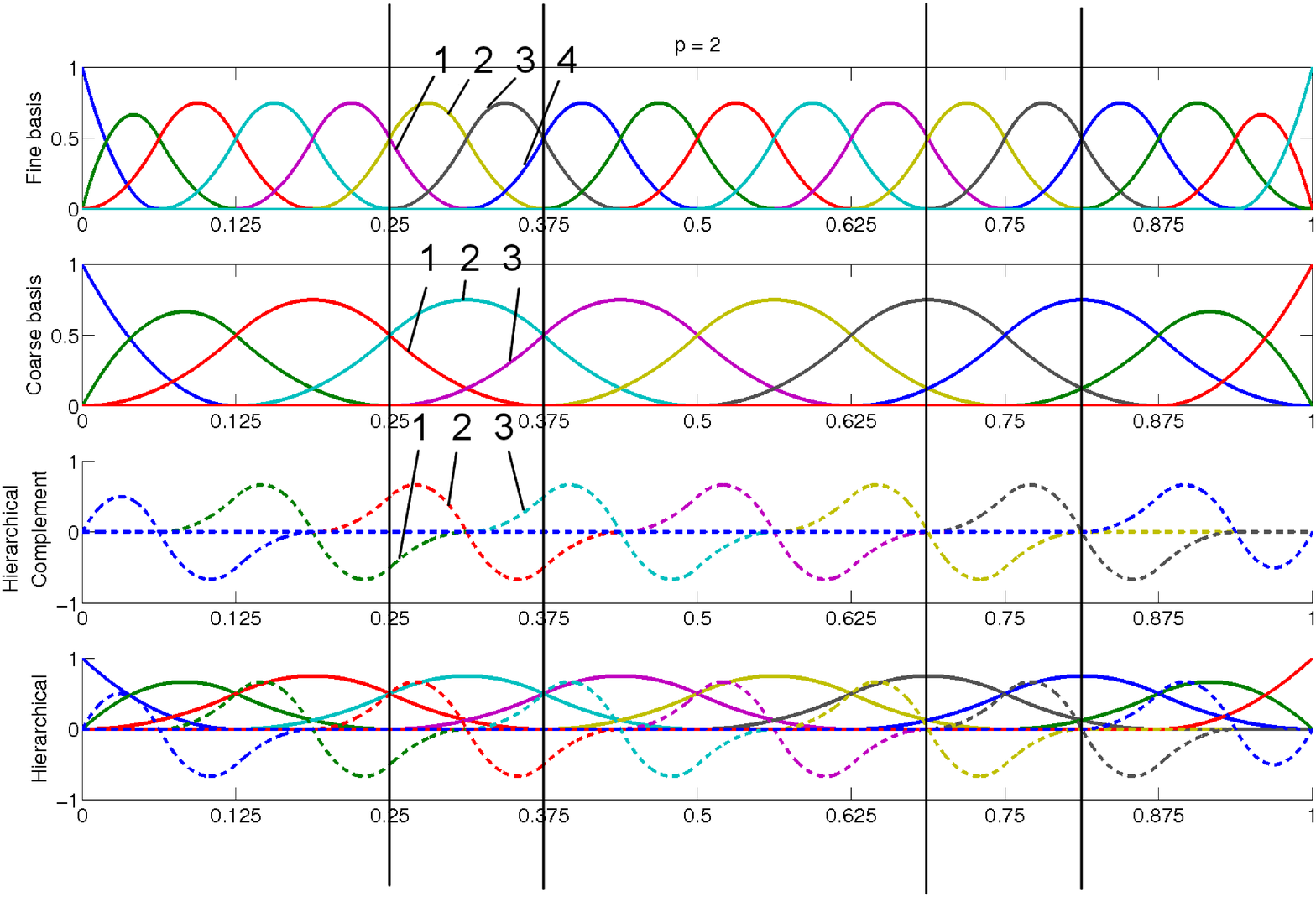}
\end{figure}

\begin{table}[p]
\begin{center}
\caption{$\gamma^2$ with first choice of $T^{p,r}_{{k}}$ in a square domain}
\label{tab:gamma_Sparse_Square}
\begin{tabular}{|c|c|c|c|c|c|}\hline
$1/h$ & 8 & 16 & 32 & 64 & 128\\\hline
$p=2,(C^{p-1})$ & 0.18 & 0.19 & 0.19 & 0.19 & 0.19\\\hline
$p=2,(C^0)$  & 0.29 & 0.32 & 0.32 & 0.32 & 0.32\\\hline
$p=3,(C^{p-1})$ & 0.36 & 0.30 & 0.30 & 0.30 & 0.30\\\hline
$p=3,(C^0)$  & 0.56 & 0.57 & 0.57 & 0.58 & 0.58\\\hline
$p=4,(C^{p-1})$ & 0.53 & 0.53 & 0.51 & 0.51 & 0.51\\\hline
$p=4,(C^0)$  & 0.78 & 0.79 & 0.79 & 0.79 & 0.79\\\hline
\end{tabular}
\end{center}
\end{table}

\begin{table}[p]
\begin{center}
\caption{$\gamma^2$ with second choice of $T^{p,r}_{{k}}$ in a square domain}
\label{tab:gamma_Dense_Square}
\begin{tabular}{|c|c|c|c|c|c|}\hline
$1/h$ & 8 & 16 & 32 & 64 & 128\\\hline
$p=2,(C^{p-1})$ & 0.09&  0.08& 0.08 & 0.08 & 0.08\\\hline
$p=2,(C^0)$ & 0.27 & 0.28 &  0.29&  0.29 & 0.29 \\\hline
$p=3,(C^{p-1})$ & 0.19 & 0.18 & 0.18 & 0.18 & 0.18\\\hline
$p=3,(C^0)$ & 0.32 & 0.33 & 0.34 & 0.34 & 0.34\\\hline
$p=4,(C^{p-1})$ & 0.53 & 0.53 & 0.51 & 0.51 & 0.51\\\hline
$p=4,(C^0)$ & 0.41 & 0.42 & 0.42 & 0.42 & 0.42\\\hline
\end{tabular}
\end{center}
\end{table}

\begin{table}
\begin{center}
\caption{$\gamma^2$ with first choice of $T^{p,r}_{{k}}$ in a quarter annulus domain}
\label{tab:gamma_Sparse_Ring}
\begin{tabular}{|c|c|c|c|c|c|}\hline
$1/h$ & 8 & 16 & 32 & 64 & 128\\\hline
$p=2,(C^{p-1})$ & 0.28 & 0.29 & 0.30 & 0.30 & 0.30\\\hline
$p=2,(C^0)$   & 0.52 & 0.56 & 0.57 & 0.58 & 0.58\\\hline
$p=3,(C^{p-1})$ & 0.44 & 0.38 & 0.38 & 0.38 & 0.38\\\hline
$p=3,(C^0)$   & 0.65 & 0.67 & 0.68 & 0.68 & 0.68\\\hline
$p=4,(C^{p-1})$ & 0.60 & 0.60 & 0.58 & 0.58 & 0.58\\\hline
$p=4,(C^0)$   & 0.85 & 0.85 & 0.85 & 0.86 & 0.86\\\hline
\end{tabular}
\end{center}
\end{table}

\begin{table}
\begin{center}
\caption{$\gamma^2$ with second choice of $T^{p,r}_{{k}}$ in a quarter annulus domain}
\label{tab:gamma_Dense_Ring}
\begin{tabular}{|c|c|c|c|c|c|}\hline
$1/h$ & 8 & 16 & 32 & 64 & 128\\\hline
$p=2,(C^{p-1})$ & 0.12 & 0.11 & 0.11 & 0.12 & 0.12\\\hline
$p=2,(C^0)$ & 0.44 & 0.47 & 0.48 & 0.49 & 0.49\\\hline
$p=3,(C^{p-1})$ & 0.29 & 0.27 & 0.26 & 0.25 & 0.25\\\hline
$p=3,(C^0)$ & 0.52 & 0.56 & 0.58 & 0.58 & 0.58\\\hline
$p=4,(C^{p-1})$ & 0.60 & 0.60 & 0.58 & 0.58 & 0.57\\\hline
$p=4,(C^0)$ & 0.53 & 0.55 & 0.57 & 0.57 & 0.57\\\hline
\end{tabular}
\end{center}
\end{table}

\begin{table}
\begin{center}
\caption{$\kappa(\hat{A}_{11})$ with first choice of $T^{p,r}_{{k}}$ in a square domain}
\label{tab:kappaA11_Sparse_Square}
\begin{tabular}{|c|c|c|c|c|c|}\hline
$1/h$ & 8 & 16 & 32 & 64 & 128\\\hline
$p=2,(C^{p-1})$ & 6.5 & 6.5 & 6.4 & 6.4 & 6.5\\\hline
$p=2,(C^0)$ & 15.9 & 17.0 & 17.3 & 17.3 & 17.4\\\hline
$p=3,(C^{p-1})$ & 24.6 & 27.3 & 28.4 & 28.8 & 28.9\\\hline
$p=3,(C^0)$ & 49.8 & 51.4 & 51.9 & 52.0 & 52.0\\\hline
$p=4,(C^{p-1})$ & 101.3 & 107.8 & 108.6 & 110.3 & 110.9\\\hline
$p=4,(C^0)$ & 322.5 & 333.7 & 336.6 & 336.6 & 337.5\\\hline
\end{tabular}
\end{center}
\end{table}

In the framework of isogeometric analysis, the local analysis of the CBS constant for $C^0$-continuous basis functions can be done as in the finite element analysis. However, for $C^{p-1}$-continuous basis functions, it is not straightforward. The extended support of B-splines (NURBS) in general creates
a dimension mismatch between the fine space and its corresponding hierarchical space on the macro element level. That is, for a given macro element the number of basis functions in the fine space is not
identical with the number of basis functions in its hierarchical space. This problem is illustrated and explained for $p=2$ in Fig. \ref{fig:HSpace} and Fig. \ref{fig:Dimension}. A local analysis in case of $C^{p-1}$ continuous basis functions therefore requires further investigations.

In Tables~\ref{tab:gamma_Sparse_Square}-\ref{tab:gamma_Dense_Ring}, we provide the global $\gamma$,
and in Tables~\ref{tab:kappaA11_Sparse_Square}-\ref{tab:kappaA11_Dense_Ring} the condition number of $\hat{A}_{11}$ is presented.
In the AMLI framework these are two decisive quantities for assessing the quality of the hierarchical
two-level splitting. The results presented in Tables~\ref{tab:kappaA11_Sparse_Square}-\ref{tab:kappaA11_Dense_Ring} show that the condition number of $\hat{A}_{11}$ block is independent of $h$. In Table~\ref{tab:kappaA11_Dense_Square} and Table~\ref{tab:kappaA11_Dense_Ring}, the entries marked by $*$ represent the cases where the results could not be obtained due to limitation on computational resources.

Note that for $p=4$ and $C^0$-continuity, the value of $\gamma^2$ in Table~\ref{tab:gamma_Sparse_Square}
and Table~\ref{tab:gamma_Sparse_Ring} is not less than $3/4$, however, as can be seen from
the numerical tests presented in the next section, the $W$-cycle still resulted in a uniform
preconditioner.

\section{Numerical results }
\label{sec:NumRes}

To test the performance of the proposed AMLI methods for IGA, we consider the following test problems, whose discretizations are performed using the Matlab toolbox GeoPDEs \cite{FalcoRV-11,GeoPDEs}.
\begin{example}
\label{exmp:Ex1}
Let $\Omega = (0,1)^{2}$. Together with $\mathcal{A} = I$, and  Dirichlet boundary conditions, the right hand side function $f$ is chosen such that the analytical solution of the problem is given by $u = e^x \sin (y)$.
\end{example}
\begin{example}
\label{exmp:Ex2}
The domain is chosen as a quarter annulus in the first Cartesian quadrant with inner radius $1$ and outer radius $2$. Together with $\mathcal{A} = I$, and homogeneous Dirichlet boundary conditions, the right hand side function $f$ is chosen such that the analytic solution is given by $u = -xy^2(x^{2} + y^{2}-1) (x^{2} + y^{2}-4)$, see  \cite{FalcoRV-11,GeoPDEs}.
\end{example}
\begin{example}
\label{exmp:Ex3}
The domain is chosen as a quarter of a thick ring. Together with $\mathcal{A} = I$, and  Dirichlet boundary conditions, the right hand side function $f$ is chosen such that the analytical solution of the problem is given by $u = e^x \sin (x y) cos ( z)$.
\end{example}

\begin{table}
\begin{center}
\caption{$\kappa(\hat{A}_{11})$ with second choice of $T^{p,r}_{{k}}$ in a square domain}
\label{tab:kappaA11_Dense_Square}
\begin{tabular}{|c|c|c|c|c|c|}\hline
$1/h$ & 8 & 16 & 32 & 64 & 128\\\hline
$p=2,(C^{p-1})$ & 14.2 & 15.0 & 15.2 & 15.3 & 15.3\\\hline
$p=2,(C^{0})$ & 20.2 & 28.8 & 33.6 & 34.9 & 35.1\\\hline
$p=3,(C^{p-1})$ & 31.6 & 42.1 & 43.4 & 43.6 & 43.8\\\hline
$p=3,(C^0)$  & 306.2 & 321.1& 325.5 & 326.5 & *\\\hline
$p=4,(C^{p-1})$ & 101.3 & 107.7 & 108.6 & 110.3 & 110.9\\\hline
$p=4,(C^0)$  & 1392.2 & 1437.1  & 1449.1 & 1452.1 & *\\\hline
\end{tabular}
\end{center}
\end{table}

\begin{table}
\begin{center}
\caption{$\kappa(\hat{A}_{11})$ with first choice of $T^{p,r}_{{k}}$ in a quarter annulus domain}
\label{tab:kappaA11_Sparse_Ring}
\begin{tabular}{|c|c|c|c|c|c|}\hline
$1/h$ & 8 & 16 & 32 & 64 & 128\\\hline
$p=2,(C^{p-1})$ & 20.1 & 22.5 & 23.6 & 24.2 & 24.5\\\hline
$p=2,(C^0)$  & 40.0 & 45.4 & 48.8 & 50.9 & 52.3\\\hline
$p=3,(C^{p-1})$ & 57.4 & 71.9 & 79.6 & 84.0 & 86.6\\\hline
$p=3,(C^0)$  & 143.1 & 155.0 & 161.6 & 165.4 & 167.4\\\hline
$p=4,(C^{p-1})$ & 220.9 & 269.8 & 298.3 & 319.0 & 331.4\\\hline
$p=4,(C^0)$  & 896.0 & 973.3 & 1007.5 & 1027.7 & 1041.6\\\hline
\end{tabular}
\end{center}
\end{table}

\begin{table}
\begin{center}
\caption{$\kappa(\hat{A}_{11})$ with second choice of $T^{p,r}_{{k}}$ in a quarter annulus domain}
\label{tab:kappaA11_Dense_Ring}
\begin{tabular}{|c|c|c|c|c|c|}\hline
$1/h$ & 8 & 16 & 32 & 64 & 128\\\hline
$p=2,(C^{p-1})$ & 43.8 & 65.4 & 81.3 & 91.4 & 98.0 \\\hline
$p=2,(C^{0})$ & 39.6 & 46.0 & 49.6 & 51.5&  52.5 \\\hline
$p=3,(C^{p-1})$ & 74.8 & 109.9 & 127.8 & 137.1 & 142.2\\\hline
$p=3,(C^{0})$ & 787.0& 870.8 & 926.7 & 965.7 & * \\\hline
$p=4,(C^{p-1})$ & 220.9 & 269.8 & 298.3 & 319.0 & 331.4\\\hline
$p=4,(C^{0})$ & 4161.5& 4561.5 & 4751.9 & 4848.1& * \\\hline
\end{tabular}
\end{center}
\end{table}

At the finest level (largest problem size), the parametric domain is divided into $n$ equal elements in each direction.
The initial guess for (iteratively) solving the linear system of equations is chosen as the zero vector. Let $r_{\mathrm{0}}$ denote the initial residual vector and $r_{\mathrm{it}}$ denote the residual vector at a given PCG/FCG iteration $n_{\mathrm{it}}$. The following stopping criteria is used
\begin{equation}
\frac{\Vert r_{\mathrm{it}} \Vert}{\Vert r_{0} \Vert} \le 10^{-8}.
\end{equation}
The average convergence factor reported in the following tables is defined as $\rho = \Big ( \dfrac{\Vert r_{\mathrm{it}} \Vert}{\Vert r_{0} \Vert}\Big )^{1/n_{\mathrm{it}}}$.
In the following tables, by $\mathrm{L1}$, $\mathrm{L2}$ and $\mathrm{N2}$ we denote the linear multiplicative AMLI cycles with $\nu=1$, $\nu=2$ and non-linear multiplicative AMLI cycle with $\nu=2$, respectively. 
By $t_c$, we represent the setup time in seconds, i.e., the time taken in the construction of transfer operators and generating the preconditioner for $\hat{A}_{11}$ block (for which we used the ILU(0) factorization, i.e. without any fill-in). The solver time (in seconds) is represented by $t_s$. All the numerical tests are performed on Intel\textsuperscript{\textregistered} Xeon\textsuperscript{\textregistered} CPU E5-1650 $@$ 3.2GHz 12 Cores and 16GB RAM. For two- and three-dimensional examples, at the coarsest level we have $h=1/4$ and $h=1/2$, respectively.  Therefore in two-dimensions, we refine the mesh upto $7$-levels of refinement and in three-dimensions upto $5$-levels of refinement. For all the test cases we take the polynomial degree $p = 2,3,4$ with $C^0$- and $C^{p-1}$-continuity. Furthermore, the transfer operator $G^{p,r}_{{k}}$ is fixed and it exactly  represents the coarse basis functions in the space of fine basis functions. 
The hierarchical complementary transfer operator $T^{p,r}_{{k}}$ are chosen in two different ways as defined in Section \ref{sec:hier_space}, see \eqref{eq:TFCp}-\eqref{eq:TSC0}.
We first consider the Example 1 and provide  $t_c$, $t_s$, $n_{it}$ and $\rho$ for $\mathrm{L1}$-, $\mathrm{L2}$-, $\mathrm{N2}$- cycles with both the choices of  $T^{p,r}_{{k}}$. Numerical results are presented in Tables \ref{tab:Ex1FCp}-\ref{tab:Ex1FC0} and Tables \ref{tab:Ex1SCp}-\ref{tab:Ex1SC0} for the first choice and the second choice of $T^{p,r}_{{k}}$, respectively .
\begin{table}[h]
{%
\begin{center}
\caption{AMLI methods for Example 1: First choice of $T^{p,r}_{{k}}$ given in \eqref{eq:TFCp} with $C^{p-1}$ regularity}
\label{tab:Ex1FCp}
\begin{tabular}{|c|c|c|c|c|c|c|c|c|c|c|}\hline
$1/h$& {$t_c$} & \mc{3}{|c|}{$t_s$} & \mc{3}{|c|}{$n_{it}$} & \mc{3}{|c|}{$\rho$}\\\hline
 &  & $\mathrm{L1}$ & $\mathrm{L2}$ & $\mathrm{N2}$ & $\mathrm{L1}$ & $\mathrm{L2}$ & $\mathrm{N2}$ & $\mathrm{L1}$ & $\mathrm{L2}$ & $\mathrm{N2}$\\\hline
\multicolumn{11}{|c|}{$p=2$} \\\hline
8 & 0.00 & 0.00 & 0.00 & 0.00 & 7 & 7 & 7 & 0.0641 & 0.0641 & 0.0622\\
16 &  0.00 & 0.00 & 0.01 & 0.01 & 8 & 7 & 7 & 0.0948 & 0.0966 & 0.0670\\
32 &  0.01 & 0.01 & 0.01 & 0.01 & 9 & 8 & 7 & 0.1108 & 0.0988 & 0.0672\\
64 &  0.04 & 0.02 & 0.03 & 0.04 & 9 & 8 & 7 & 0.1086 & 0.0901 & 0.0622\\
128 & 0.18 & 0.07 & 0.09 & 0.12 & 9 & 8 & 7 & 0.1166 & 0.0909 & 0.0624\\
256  & 0.72 & 0.25 & 0.30 & 0.41 & 9 & 8 & 7 & 0.1175 & 0.0879 & 0.0603\\
512 & 2.97 & 1.03 & 1.12 & 1.50 & 9 & 8 & 7 & 0.1276 & 0.0945 & 0.0620\\\hline
\multicolumn{11}{|c|}{$p=3$} \\\hline
8 & 0.00 & 0.00 & 0.00 & 0.00 & 8 & 8 & 8 & 0.0901 & 0.0901 & 0.0901\\
16  & 0.01 & 0.01 & 0.01 & 0.01 & 9 & 9 & 8 & 0.1111 & 0.1129 & 0.0686\\
32 & 0.02 & 0.01 & 0.01 & 0.02 & 10 & 9 & 7 & 0.1293 & 0.1043 & 0.0577\\
64 &  0.10 & 0.03 & 0.04 & 0.05 & 10 & 8 & 7 & 0.1361 & 0.0857 & 0.0551\\
128 & 0.41 & 0.12 & 0.13 & 0.18 & 10 & 8 & 7 & 0.1369 & 0.0821 & 0.0536\\
256 &  1.76 & 0.48 & 0.46 & 0.63 & 10 & 8 & 7 & 0.1348 & 0.0794 & 0.0523\\
512 & 7.50 & 1.65 & 1.77 & 2.37 & 9 & 8 & 7 & 0.1283 & 0.0771 & 0.0511\\\hline
\multicolumn{11}{|c|}{$p=4$} \\\hline
8 & 0.00 & 0.00 & 0.00 & 0.00 & 10 & 10 & 10 & 0.1139 & 0.1139 & 0.1139\\
16  & 0.01 & 0.01 & 0.01 & 0.01 & 12 & 12 & 10 & 0.1866 & 0.1882 & 0.1378\\
32 & 0.06 & 0.02 & 0.02 & 0.03 & 12 & 11 & 9 & 0.2013 & 0.1822 & 0.1100\\
64 & 0.26 & 0.07 & 0.07 & 0.10 & 12 & 10 & 9 & 0.2038 & 0.1557 & 0.1032\\
128 & 1.09 & 0.26 & 0.24 & 0.37 & 12 & 9 & 9 & 0.2028 & 0.1209 & 0.0977\\
256 & 4.57 & 0.98 & 0.88 & 1.21 & 12 & 9 & 8 & 0.1976 & 0.1182 & 0.0975\\
512 &19.05 & 3.60 & 3.44 & 4.66 & 11 & 9 & 8 & 0.1853 & 0.1146 & 0.0930\\\hline
\end{tabular}
\end{center}
}
\end{table}
\begin{table}[t!]
{%
\begin{center}
\caption{AMLI methods for Example 1:  First choice of $T^{p,r}_{{k}}$ given in \eqref{eq:TFC0} with  $C^{0}$ regularity}
\label{tab:Ex1FC0}
\begin{tabular}{|c|c|c|c|c|c|c|c|c|c|c|}\hline
$1/h$& {$t_c$} & \mc{3}{|c|}{$t_s$} & \mc{3}{|c|}{$n_{it}$} & \mc{3}{|c|}{$\rho$}\\\hline
 &  & $\mathrm{L1}$ & $\mathrm{L2}$ & $\mathrm{N2}$ & $\mathrm{L1}$ & $\mathrm{L2}$ & $\mathrm{N2}$ & $\mathrm{L1}$ & $\mathrm{L2}$ & $\mathrm{N2}$\\\hline
\multicolumn{11}{|c|}{$p=2$} \\\hline
8 & 0.00 & 0.01 & 0.01 & 0.01 & 9 & 9 & 9 & 0.1072 & 0.1072 & 0.1072\\
16 & 0.01 & 0.01 & 0.01 & 0.01 & 11 & 11 & 9 & 0.1695 & 0.1716 & 0.1102\\
32  & 0.02 & 0.03 & 0.04 & 0.04 & 13 & 11 & 9 & 02195 & 0.1738 & 0.1110\\
64  & 0.10 & 0.07 & 0.09 & 0.10 & 14 & 11 & 9 & 0.2606 & 0.1744 & 0.1109\\
128  & 0.38 & 0.30 & 0.29 & 0.34 & 16 & 11 & 9 & 0.2973 & 0.1743 & 0.1105\\
256 & 1.65 & 1.25 & 1.04 & 1.23 & 17 & 11 & 9 & 0.3288 & 0.1736 & 0.1102\\
512  & 6.93 & 5.17 & 3.84 & 4.61 & 18 & 11 & 9 & 0.3557 & 0.1730 & 0.1100\\\hline
\multicolumn{11}{|c|}{$p=3$ } \\\hline
8  & 0.01 & 0.01 & 0.01 & 0.03 & 12 & 12 & 12 & 0.1999&  0.1999& 0.1999\\
16  & 0.02 & 0.03 & 0.04 & 0.03 & 17& 17 & 12 & 0.3288 & 0.3305 & 0.2124\\
32& 0.09 & 0.09 & 0.12 & 0.10& 22 & 18& 12 & 0.4258 &0.3568 & 0.2129\\
64 &  0.37 & 0.38 & 0.41 & 0.33 &27 & 19 & 12 & 0.5014 & 0.3650 &0.2122\\
128& 1.55& 1.77 & 1.34 & 1.21 &32 & 19 & 12 & 0.5581& 0.3673 & 0.2114\\
256 &  6.73 & 7.87 & 4.88 & 4.51 &37 & 19 & 12& 0.6038 & 0.3670 & 0.2110\\
512 & 28.76 & 36.56 & 19.12 & 17.84 &42 & 19 & 12& 0.6394 & 0.3664 & 0.2108\\\hline
\multicolumn{11}{|c|}{$p=4$ } \\\hline
8      & 0.01 & 0.03 & 0.03 & 0.03 & 19 & 19 & 19 & 0.3631 & 0.3631 & 0.3631\\
16    & 0.05 & 0.07 & 0.10 & 0.09 & 25 & 26 & 19 & 0.4784 & 0.4827 & 0.3719\\
32    & 0.24 & 0.32 & 0.36 & 0.29 & 38 & 30 & 19 & 0.6087 & 0.5337 & 0.3720\\
64    & 1.07 & 1.62 & 1.29 & 1.05 & 52 & 32 & 19 & 0.6982 & 0.5585 & 0.3719\\
128  & 4.50 & 8.04 & 4.90 & 3.89 & 67 & 34 & 19 & 0.7585 & 0.5766 & 0.3709\\
256  & 18.73 & 40.11 & 19.41 & 15.25 & 85 & 35 & 19 & 0.8038 & 0.5827 & 0.3703\\
512  & 76.22& 190.29& 77.37& 62.24& 100\footnotemark&35 & 19& 0.8379 & 0.5878 & 0.3700\\\hline
\end{tabular}
\end{center}
}
\end{table}
\begin{table}[h]
\begin{center}
\caption{AMLI methods for Example 1: Second choice of $T^{p,r}_{{k}}$ given in \eqref{eq:TSCp} with  $C^{p-1}$ regularity}
\label{tab:Ex1SCp}
\begin{tabular}{|c|c|c|c|c|c|c|c|c|c|c|}\hline
$1/h$& {$t_c$} & \mc{3}{|c|}{$t_s$} & \mc{3}{|c|}{$n_{it}$} & \mc{3}{|c|}{$\rho$}\\\hline
 &  & $\mathrm{L1}$ & $\mathrm{L2}$ & $\mathrm{N2}$ & $\mathrm{L1}$ & $\mathrm{L2}$ & $\mathrm{N2}$ & $\mathrm{L1}$ & $\mathrm{L2}$ & $\mathrm{N2}$\\\hline
\multicolumn{11}{|c|}{$p=2$ } \\\hline
8 &  0.08 & 0.02 & 0.42 & 0.52 & 5 & 5 & 5 & 0.0227 & 0.0227 & 0.0227\\
16 &  0.00 & 0.01 & 0.01 & 0.01 & 6 & 6 & 5 & 0.0304 & 0.0326 & 0.0217\\
32 &  0.02 & 0.01 & 0.01 & 0.01 & 6 & 6 & 5 & 0.0316 & 0.0311 & 0.0226\\
64  & 0.07 & 0.02 & 0.05 & 0.05 & 6 & 6 & 5 & 0.0303 & 0.0300 & 0.0224\\
128 &  0.30 & 0.06 & 0.08 & 0.10 & 6 & 6 & 5 & 0.0314 & 0.0310 & 0.0234\\
256 & 1.21 & 0.22 & 0.30 & 0.39 & 6 & 6 & 5 & 0.0301 & 0.0296 & 0.0226\\
512 &  5.18 & 0.88 & 1.05 & 1.62 & 6 & 6 & 6 & 0.0326 & 0.0321 & 0.0269\\\hline
\multicolumn{11}{|c|}{$p=3$} \\\hline
8 & 0.00 & 0.02 & 0.00 & 0.00 & 7 & 7 & 7 & 0.0443 & 0.0443 & 0.0443\\
16 &  0.01 & 0.00 & 0.00 & 0.01 & 7 & 7 & 6 & 0.0560 & 0.0569 & 0.0365\\
32 & 0.04 & 0.01 & 0.01 & 0.02 & 7 & 7 & 6 & 0.0576 & 0.0494 & 0.0319\\
64 &  0.18 & 0.04 & 0.04 & 0.05 & 7 & 6 & 5 & 0.0569 & 0.0377 & 0.0216\\
128 &  0.77 & 0.12 & 0.13 & 0.17 & 7 & 6 & 5 & 0.0542 & 0.0343 & 0.0204\\
256 & 3.32 & 0.47 & 0.49 & 0.64 & 7 & 6 & 5 & 0.0502 & 0.0326 & 0.0195\\
512 &  13.99 & 1.60 & 1.89 & 2.44 & 6 & 6 & 5 & 0.0446 & 0.0311 & 0.0186\\\hline
\end{tabular}
\end{center}
\end{table}
\begin{table}[h]
\begin{center}
\caption{AMLI methods for Example 1: Second choice of $T^{p,r}_{{k}}$ given in \eqref{eq:TSC0} with  $C^0$ regularity}
\label{tab:Ex1SC0}
\begin{tabular}{|c|c|c|c|c|c|c|c|c|c|c|}\hline
$1/h$& {$t_c$} & \mc{3}{|c|}{$t_s$} & \mc{3}{|c|}{$n_{it}$} & \mc{3}{|c|}{$\rho$}\\\hline
 &  & $\mathrm{L1}$ & $\mathrm{L2}$ & $\mathrm{N2}$ & $\mathrm{L1}$ & $\mathrm{L2}$ & $\mathrm{N2}$ & $\mathrm{L1}$ & $\mathrm{L2}$ & $\mathrm{N2}$\\\hline
\multicolumn{11}{|c|}{$p=2$} \\\hline
8 & 0.00 & 0.00 & 0.00 & 0.00 & 8 & 8 & 8 & 0.0901 & 0.0901 & 0.0901\\
16 & 0.01 & 0.01 & 0.01 & 0.01 & 9 & 9 & 8 & 0.1173 & 0.1195 & 0.0918\\
32 & 0.04 & 0.02 & 0.03 & 0.04 & 10 & 9 & 8 & 0.1308 & 0.1197 & 0.0900\\
64  & 0.15 & 0.07 & 0.09 & 0.12 & 10 & 9 & 8 & 0.1430 & 0.1192 & 0.0890\\
128  & 0.62 & 0.27 & 0.32 & 0.41 & 10 & 9 & 8 & 0.1479 & 0.1191 & 0.0884\\
256  & 2.71 & 1.08 & 1.20 & 1.56 & 10 & 9 & 8 & 0.1509 & 0.1191 & 0.0880\\
512 & 11.22 & 4.33 & 4.55 & 6.00 & 10 & 9 & 8 & 0.1528 & 0.1191 & 0.0878\\\hline
\multicolumn{11}{|c|}{$p=3$} \\\hline
8 & 0.01 & 0.01 & 0.03 & 0.03 & 9 & 9 & 9 & 0.1133 & 0.1133 & 0.1133\\
16 & 0.02 & 0.02 & 0.02 & 0.02 & 11 & 11 & 9 & 0.1724 & 0.1743 & 0.1191\\
32 & 0.09 & 0.05 & 0.07 & 0.07 & 13 & 11 & 9 & 0.2216 & 0.1762 & 0.1206\\
64 & 0.39 & 0.20 & 0.22 & 0.25 & 14 & 11 & 9 & 0.2627 & 0.1777 & 0.1212\\
128 &  1.62 & 0.88 & 0.79 & 0.91 & 16 & 11 & 9 & 0.2998 & 0.1782 & 0.1215\\
256 &  7.06 & 3.70 & 2.87 & 3.42 & 17 & 11 & 9 & 0.3321 & 0.1785 & 0.1216\\
512 & 29.78 & 16.54 & 11.24 & 13.37 & 19 & 11 & 9 & 0.3630 & 0.1786 & 0.1217\\\hline
\multicolumn{11}{|c|}{$p=4$} \\\hline
8 & 0.01 & 0.02 & 0.02 & 0.02 & 10 & 10 & 10 & 0.1368 & 0.1368 & 0.1368\\
16 & 0.07 & 0.04 & 0.06 & 0.05 & 13 & 13 & 10 & 0.2199 & 0.2219 & 0.1419\\
32  & 0.33 & 0.15 & 0.18 & 0.18 & 15 & 13 & 10 & 0.2825 & 0.2254 & 0.1416\\
64  & 1.88 & 0.61 & 0.60 & 0.64 & 17 & 13 & 10 & 0.3259 & 0.2251 & 0.1412\\
128  & 5.92 & 2.51 & 2.20 & 2.41 & 18 & 13 & 10 & 0.3575 & 0.2247 & 0.1410\\
256& 25.51 & 11.39 & 8.56 & 9.56 & 20 & 13 & 10 & 0.3844 & 0.2245 & 0.1408\\
512   & 104.33 & 49.49 & 35.15 & 39.67& 21& 13 & 10& 0.4042 & 0.2244 & 0.1408\\\hline
\end{tabular}
\end{center}
\end{table}
From Tables \ref{tab:Ex1FCp}-\ref{tab:Ex1SC0} we observe the following:
\begin{itemize}

\item The number of iterations and total solution $(t_c+t_s)$ time show an  $h$-independent convergence rates for $C^{p-1}$- and $C^0$-continuity.

\item For $C^{p-1}$-continuity, the results are  $p$-independent, whereas for $C^0$-continuity, the degree $p$ has some effect on PCG/FCG iterations.

\item For $C^{p-1}$-continuity, all the AMLI cycles give optimal results, and the $V$-cycle $(\nu=1)$ is the fastest among all. This is due to a very nice bound on $\gamma$ for $C^{p-1}$-continuity. Among $\mathrm{L2}$-, $\mathrm{N2}$-cycles, the latter has smaller iteration numbers. Therefore, in the remaining numerical computations we consider linear AMLI cycle with $\nu=1$ and nonlinear AMLI cycle with $\nu=2$ for $C^{p-1}$ continuous basis functions.

\item For $C^0$-continuity,  $V$-cycle $(\nu=1)$ is not an optimal order method, an observation similar to standard FEM. However, for $C^0$-continuity, both the $\nu=2$ cycle methods (linear and nonlinear) exhibit optimal order behavior, and nonlinear AMLI gives overall better results. Therefore, we consider only nonlinear AMLI cycle with $\nu=2$ for $C^0$ continuous basis functions in remaining numerical results. 

\item For $p=4$ with $C^{p-1}$-continuity, we could not obtain better $\gamma$ with the second choice of $T^{p,r}_{{k}}$ as compared to the first choice. Therefore, in Table \ref{tab:Ex1SCp}, the numerical results are presented only for $p=2,3$ with second choice of $T^{p,r}_{{k}}$. Numerical results for $p=4$ may be improved by choosing different operators, which demands further investigation.

 \item For $C^{p-1}$-continuity, though the number of iterations are less for second choice of $T^{p,r}_{{k}}$, the overall time $(t_c+t_s)$ is more than the first choice of $T^{p,r}_{{k}}$. This happens due to comparatively less sparse structure of second choice $T^{p,r}_{{k}}$, which results in more construction time $t_c$. Therefore, in the remaining numerical tests we consider only the first choice of  $T^{p,r}_{{k}}$ for $C^{p-1}$ continuous basis functions.

\item For $C^0$-continuity, we get mixed results from both the choices of $T^{p,r}_{{k}}$. This is due to the fact that there is not much difference in number of nonzero entries in each row of $T^{p,r}_{{k}}$ for two different choices. Therefore, numerical results are provided for both the choices of $T^{p,r}_{{k}}$ for $C^0$ continuous basis functions. 

\end{itemize}

We now consider Example 2 with curved boundary. The geometry for this example is represented by NURBS basis functions of order $1$ in the radial direction and of order $2$ in the angular direction, see \cite{FalcoRV-11}. Numerical results are provided for $C^{p-1}$-continuity with first choice of $T^{p,r}_{{k}}$ in Table \ref{tab:Ex2FCp}, and for $C^0$-continuity with both the choices of $T^{p,r}_{{k}}$ in Table \ref{tab:Ex2FCSC0}. All the results are qualitatively similar to that of Example $1$ with square domain.
\begin{table}[h]
\caption{AMLI methods for Example 2: First choice of $T^{p,r}_{{k}}$ given in \eqref{eq:TFCp} with  $C^{p-1}$ regularity}
\label{tab:Ex2FCp}
\begin{center}
\begin{tabular}{|c|c|c|c|c|c|c|c|}\hline
$1/h$ & $t_c$ & \mc{2}{|c|}{$t_s$} & \mc{2}{|c|}{$n_{it}$} & \mc{2}{|c|}{$\rho$}\\\hline
&  & $\mathrm{L1}$ &  $\mathrm{N2}$ & $\mathrm{L1}$ & $\mathrm{N2}$ &$\mathrm{L1}$ & $\mathrm{N2}$\\\hline
\multicolumn{8}{|c|}{$p=2$} \\\hline
8 & 0.02 & 0.02 & 0.01 & 8 & 8 & 0.0802 & 0.0802\\
16  & 0.00 & 0.01 & 0.01 & 9 & 8 & 0.1201 & 0.0839\\
32 & 0.01 & 0.01 & 0.01 & 10 & 7 & 0.1499 & 0.0658\\
64 & 0.05 & 0.02 & 0.03 & 11 & 6 & 0.1838 & 0.0453\\
128 & 0.17 & 0.09 & 0.10 & 12 & 6 & 0.2048 & 0.0351\\
256 & 0.72 & 0.38 & 0.30 & 13 & 5 & 0.2211 & 0.0226\\
512  & 2.93 & 1.53 & 1.07 & 13 & 5 & 0.2374 & 0.0194\\\hline
\multicolumn{8}{|c|}{$p=3$} \\\hline
8 & 0.00 & 0.00 & 0.00 & 9 & 9 & 0.1201 & 0.1201\\
16 & 0.01 & 0.01 & 0.01 & 10 & 9 & 0.1560 & 0.1148\\
32 & 0.02 & 0.01 & 0.02 & 12 & 8 & 0.1839 & 0.0988\\
64 & 0.10 & 0.04 & 0.06 & 13 & 8 & 0.2104 & 0.0900\\
128 & 0.41 & 0.16 & 0.20 & 13 & 8 & 0.2363 & 0.0858\\
256 & 1.76 & 0.66 & 0.72 & 14 & 8 & 0.2514 & 0.0828\\
512 & 7.45 & 2.56 & 2.35 & 14 & 7 & 0.2644 & 0.0706\\\hline
\multicolumn{8}{|c|}{$p=4$} \\\hline
8 &  0.03 & 0.01 & 0.00 & 11 & 11 & 0.1686 & 0.1686\\
16 &  0.01 & 0.01 & 0.01 & 12 & 11 & 0.2073 & 0.1665\\
32 &  0.05 & 0.02 & 0.03 & 13 & 9 & 0.2419 & 0.1248\\
64 &  0.24 & 0.11 & 0.10 & 14 & 9 & 0.2549  & 0.1054\\
128 &  1.07 & 0.32 & 0.43 & 15  & 8 & 0.2688  & 0.0884\\
256 &  4.47 & 1.23  & 1.09 & 15 & 7 & 0.2924  & 0.0648\\
512 &  18.79 & 5.30 & 4.13 & 16 & 7 & 0.3061 & 0.0534\\\hline
\end{tabular}
\end{center}
\end{table}
\begin{table}[h]
\caption{AMLI methods for Example 2:   $C^{0}$ regularity}
\label{tab:Ex2FCSC0}
\begin{center}
\begin{tabular}{|c|c|c|c|c|}\hline
\multicolumn{5}{|c|}{with first choice of $T^{p,r}_{{k}}$ given in \eqref{eq:TFC0}}\\\hline
$1/h$ & $t_c$ & \mc{1}{|c|}{$t_s$} & \mc{1}{|c|}{$n_{it}$} & \mc{1}{|c|}{$\rho$}\\\hline
&  &   $\mathrm{N2}$ &  $\mathrm{N2}$ & $\mathrm{N2}$\\\hline
\multicolumn{5}{|c|}{$p=2$} \\\hline
8 & 0.00 & 0.01 & 11 & 0.1744\\
16 & 0.01 & 0.02 & 11 & 0.1820\\
32 & 0.02 & 0.05 & 11 & 0.1791\\
64 & 0.09 & 0.13 & 11 & 0.1752\\
128 & 0.40 & 0.43 & 11 & 0.1730\\
256 & 1.72 & 1.52 & 11 & 0.1717\\
512 & 7.36 & 5.61 & 11 & 0.1704\\\hline
\multicolumn{5}{|c|}{$p=3$} \\\hline
8 & 0.00 & 0.01 & 13 & 0.2237\\
16 & 0.02 & 0.04 & 14 & 0.2507\\
32 & 0.08 & 0.11 & 14 & 0.2584\\
64 & 0.34 & 0.39 & 14 & 0.2632\\
128 & 1.49 & 1.43 & 14 & 0.2649\\
256 & 6.35 & 5.37 & 14 & 0.2648\\
512 & 27.51 & 20.83 & 14 & 0.2638\\\hline
\multicolumn{5}{|c|}{$p=4$} \\\hline
8 & 0.01 & 0.03 & 22 & 0.4319\\
16 & 0.05 & 0.11 & 24 & 0.4516\\
32 & 0.22 & 0.38 & 24 & 0.4563\\
64 & 0.92 & 1.34 & 24 & 0.4591\\
128 & 4.21 & 5.03 & 24 & 0.4609\\
256 & 18.28 & 19.79 & 24 & 0.4639\\
512 & 76.62 & 81.78 & 25 & 0.4644\\\hline
\end{tabular}
\quad
\begin{tabular}{|c|c|c|c|c|}\hline
\multicolumn{5}{|c|}{with second choice of $T^{p,r}_{{k}}$ given in \eqref{eq:TSC0}}\\\hline
$1/h$ & $t_c$ & \mc{1}{|c|}{$t_s$} & \mc{1}{|c|}{$n_{it}$} & \mc{1}{|c|}{$\rho$}\\\hline
&  &   $\mathrm{N2}$ &  $\mathrm{N2}$ & $\mathrm{N2}$\\\hline
\multicolumn{5}{|c|}{$p=2$} \\\hline
8 & 0.00 & 0.01 & 10 & 0.1445\\
16 &  0.01 & 0.02 & 10 & 0.1510\\
32 & 0.04 & 0.04  & 10 & 0.1478\\
64 &  0.15  & 0.15  & 10 & 0.1463\\
128 & 0.62 & 0.52  & 10 & 0.1437\\
256 &  2.65 & 1.93 & 10 &  0.1419\\
512 & 11.05 & 7.72 & 10 & 0.1401\\\hline
\multicolumn{5}{|c|}{$p=3$} \\\hline
8 &  0.01 & 0.01 & 11 & 0.1647\\
16 &0.02 & 0.03 & 11 & 0.1780\\
32 &  0.09 & 0.09 & 11 & 0.1845\\
64 &  0.39 & 0.33 & 12 & 0.1883\\
128 &  1.63 & 1.21 & 12 & 0.1922\\
256 &  6.98 & 4.52 & 12 & 0.1938\\
512 &  28.76& 17.94 & 12 & 0.1940\\\hline
\multicolumn{5}{|c|}{$p=4$} \\\hline
8 & 0.01 & 0.02 & 11 & 0.1660\\
16 &  0.07 & 0.05 & 11 & 0.1758\\
32 & 0.32 & 0.19 & 11 & 0.1789\\
64 &  1.39 & 0.70 & 11 & 0.1785\\
128  & 5.99 & 2.64 & 11 & 0.1774\\
256  & 25.31 & 10.49 & 11 & 0.1765\\
512  & 99.22 & 43.15 & 11 & 0.1757\\\hline
\end{tabular}
\end{center}
\end{table}

Finally, we consider three-dimensional problem as stated in Example 3. 
The numerical results are presented in Tables \ref{tab:Ex3FCp}-\ref{tab:Ex3FCSC0}. Due to the limitation of computer resources numerical results in three-dimensions are provided only upto $h=1/32$. In Table \ref{tab:Ex3FCp}, linear AMLI cycle with $\nu=1$, and nonlinear AMLI cycle with $\nu=2$ are given for $C^{p-1}$ continuity with first choice of  $T^{p,r}_{{k}}$. The results exhibit optimal order for both the solvers. The increased number of iterations (as compared to two-dimensional examples) can be attributed to the smaller angle between coarse space and its complementary space. 
For $C^0$-continuity the numerical results with both the choices of $T^{p,r}_{{k}}$ are given in Table \ref{tab:Ex3FCSC0}. 
The first choice of $T^{p,r}_{{k}}$, however, does not result in an optimal order method. The optimality is restored with $\nu=3$, which are presented in the column with $N3$. The second choice, though expensive, gives optimal order method for second order stabilization $(\nu=2)$.
In Tables \ref{tab:Ex3FCp}-\ref{tab:Ex3FCSC0}, The entries marked by $*$ represent the cases where the computations are performed on a  machine with larger memory but shared with other users, therefore timings are not provided for these cases.  

We note that for two-dimensional problems, the total time of the solvers also exhibit optimal complexity, however, for three-dimensional problem the increase in the total time $(t_c+t_s)$ for successive refinement is more than the factor of increase in number of unknowns. This is particularly due to two reasons, the construction of operators $G^{p,r}_{{k}}$ and $T^{p,r}_{{k}}$ by tensor product of matrices for one-dimensional operators (see Remark \ref{rem:TensorOperator}), and expensive preconditioner for $\hat{A}_{11}$ (ILU(0)). In our future study on local analysis, we also intend to construct these operators for two- and three-dimensional problems without tensor product, and devise efficient and cheaper preconditioner for $\hat{A}_{11}$.
\begin{table}[t!]
\caption{AMLI methods for Example 3: First choice of $T^{p,r}_{{k}}$ given in \eqref{eq:TFCp} with  $C^{p-1}$ regularity}
\label{tab:Ex3FCp}
\begin{center}
\begin{tabular}{|c|c|c|c|c|c|c|c|}\hline
$1/h$ & $t_c$ & \mc{2}{|c|}{$t_s$} & \mc{2}{|c|}{$n_{it}$} & \mc{2}{|c|}{$\rho$}\\\hline
&  & $\mathrm{L1}$ &  $\mathrm{N2}$ & $\mathrm{L1}$ & $\mathrm{N2}$ &$\mathrm{L1}$ & $\mathrm{N2}$\\\hline
\multicolumn{8}{|c|}{$p=2$} \\\hline
4 &  0.00 & 0.00 & 0.00 & 8 & 8 & 0.0899 & 0.0899\\
8  & 0.04 & 0.01 & 0.01 & 12 & 10 & 0.1913 & 0.1438\\
16  & 0.60 & 0.10 & 0.10 & 13 & 10 & 0.2400 & 0.1484\\
32  & 7.18 & 1.09 & 0.89 & 15 & 10 & 0.2694 & 0.1346\\
64  & * & * & * & 15 & 9 & 0.2830 & 0.1168\\\hline
\multicolumn{8}{|c|}{$p=3$} \\\hline
4  & 0.00 & 0.00 & 0.00 & 10 & 10 & 0.1415 & 0.1415\\
8 & 0.15 & 0.02 & 0.03 & 14 & 13 & 0.2492 & 0.2304\\
16 & 2.84 & 0.27 & 0.24 & 15 & 11 & 0.2923 & 0.1862\\
32 & 35.61 & 2.79 & 2.21 & 17 & 11 & 0.3215 & 0.1762\\
64  & * & *& *& 17 & 11 & 0.3349 & 0.1738\\\hline
\multicolumn{8}{|c|}{$p=4$} \\\hline
4 & 0.01 & 0.01 & 0.01 & 10 & 10 & 0.1443 & 0.1443\\
8 & 0.52 & 0.06 & 0.07 & 16 & 16 & 0.3027 & 0.3040\\
16 & 14.81 & 0.82 & 0.85 & 20 & 17 & 0.3900 & 0.3324\\
32  & 213.74 & 8.82 & 7.55 & 21 & 15 & 0.4067 & 0.2927\\
64  & * & * &* & 21 & 14 & 0.4042 & 0.2546\\\hline
\end{tabular}
\end{center}
\end{table}

\footnotetext{did not converge upto desired accuracy.}

\begin{table}[h]
\caption{AMLI methods for Example 3: with  $C^{0}$ regularity}
\label{tab:Ex3FCSC0}
\begin{center}
\begin{tabular}{|c|c|c|c|c|}\hline
\multicolumn{5}{|c|}{with first choice of $T^{p,r}_{{k}}$ given in \eqref{eq:TFC0}}\\\hline
$1/h$ & $t_c$ & \mc{1}{|c|}{$t_s$} & \mc{1}{|c|}{$n_{it}$} & \mc{1}{|c|}{$\rho$}\\\hline
&  &   $\mathrm{N2}$ &  $\mathrm{N2} (\mathrm{N3})$ & $\mathrm{N2}$\\\hline
\multicolumn{5}{|c|}{$p=2$} \\\hline
4 & 0.01 & 0.01 & 12 (12) & 0.2124\\
8 & 0.11 & 0.05 & 15 (15) & 0.2904\\
16 & 1.25 & 0.52 & 16 (15) & 0.2996\\
32 & 12.06 & 4.51 & 16 (15)& 0.3022\\\hline
\multicolumn{5}{|c|}{$p=3$} \\\hline
4 & 0.07 & 0.04 & 18 (18) & 0.3527\\
8 & 1.09 & 0.50 & 23 (22) & 0.4408\\
16 & 12.23 & 5.13 & 26 (23) & 0.4919\\
32 & 114.77 & 49.48 & 28 (23) & 0.5164\\\hline
\multicolumn{5}{|c|}{$p=4$} \\\hline
4 & 0.39 & 0.45 & 48 (48) & 0.6770\\
8 & 5.84 & 4.36 & 54 (50) & 0.7081\\
16 & 64.09 & 47.27 & 64 (51) & 0.7497\\
32 & * & * & 73 (51) & 0.7764\\\hline
\end{tabular}
\quad
\begin{tabular}{|c|c|c|c|c|}\hline
\multicolumn{5}{|c|}{with second choice of $T^{p,r}_{{k}}$ given in \eqref{eq:TSC0}}\\\hline
$1/h$ & $t_c$ & \mc{1}{|c|}{$t_s$} & \mc{1}{|c|}{$n_{it}$} & \mc{1}{|c|}{$\rho$}\\\hline
&  &   $\mathrm{N2}$ &  $\mathrm{N2}$ & $\mathrm{N2}$\\\hline
\multicolumn{5}{|c|}{$p=2$} \\\hline
4 &  0.37 & 0.31 & 11 & 0.1753\\
8  & 0.32 & 0.13 & 13 & 0.2212\\
16  & 4.29 & 0.76 & 13 & 0.2250\\
32& 33.13 & 7.44 & 13 & 0.2261\\\hline
\multicolumn{5}{|c|}{$p=3$} \\\hline
4  & 0.09 & 0.03 & 14 & 0.2663\\
8 & 1.42 & 0.34 & 16 & 0.3092\\
16 & 15.72 & 3.30 & 17 & 0.3342\\
32& 123.05 & 32.24 & 18 & 0.3415\\\hline
\multicolumn{5}{|c|}{$p=4$} \\\hline
4 &  0.98 & 0.23 & 16 & 0.2987\\
8 &  13.39 & 1.84 & 18 & 0.3465\\
16 &  144.03 & 17.32 & 18 & 0.3560\\
32& * & * & 18 & 0.3577 \\\hline
\end{tabular}
\end{center}
\end{table}

\section{Conclusions}
\label{sec:Conc}

We have presented AMLI methods for the linear system arising from the isogeometric discretization of the scalar second order elliptic problems. We summarize the main contribution of this paper as follows.
\begin{enumerate}

 \item  We provide the explicit representation of B-splines as a function of mesh size $h$ on a unit interval with uniform refinement. The explicit representation is given for $C^0$ and $C^{p-1}$ continuous basis functions of polynomial degree $p=2,3,4$, the most widely used cases in engineering applications. Explicit form of B-splines is  important from computational point of view, as well as in forming the inter-grid transfer operators.
 
\item The construction of B-spline basis functions at coarse level from the linear combination of fine basis functions is provided. For $p=2,3,4,$ and with $C^0$ and $C^{p-1}$ continuities, these transfer operators (from fine level to coarse level) are given in matrix form for a multilevel mesh. These operators can also be used to generate restriction operators in multigrid methods.

\item The transfer operators are also provided for NURBS basis functions. The formulation of NURBS operators is given in terms of B-spline operators and weights.   
 
\item The construction of hierarchical spaces for B-splines (NURBS) is presented. Hierarchical spaces are constructed as direct sum of coarse spaces and corresponding hierarchical complementary spaces. We have presented matrix form of these operators. As the choice of hierarchical complementary spaces is not unique, we have provided two different choices of these operators for each of $C^0$- and $C^{p-1}$-continuity of basis functions. 

\item For a given polynomial degree $p$, AMLI cycles are of optimal complexity with respect to the mesh refinement. Also, for a given mesh size $h$, AMLI cycles are (almost) $p$-independent. We provided numerical results for a square domain, quarter annulus (ring), and quarter thick ring. The iteration counts, convergence factor, and timings are given for AMLI linear $V$-, $W$- and nonlinear $W$-cycles. Note that, for $C^{p-1}$-continuity the linear $V$-cycle also exhibits optimal convergence (due to very nice space splitting, which is normally not found in standard FEM). The linear and nonlinear AMLI $W$-cycle is optimal for all cases except for a particular case of degree $p=4$ with $C^0$-continuity in three-dimensional problem with first choice of $T^{p,r}_k$.  For this case, the number of iterations are provided with $\nu=3$ cycle, which is optimal. The numerical results are complete for $p=2,3,4,$ with $C^{p-1}$ and $C^0$ continuous basis functions.

\end{enumerate}

Despite that the condition number of the stiffness matrix grows very rapidly with the polynomial degree, these excellent results exhibit the strength and flexibility of AMLI methods. Nevertheless, the rigorous local analysis of the CBS constant $\gamma$, particularly due to the overlapped support of B-splines,  is not a straight forward task, and is still an open problem. We intend to address this issue  in our future work.

\section*{Acknowledgement}

First two authors were partially supported by the Austrian Sciences Fund (Project \textbf{P21516-N18}). 
%

%


\end{document}